\documentclass{amsart}

\usepackage{amssymb,amscd,amsthm,latexsym}
\usepackage[all]{xy}

\def\Pt{ \mathrm{Pt} }

\def\cleft{\hbox{[\kern-.16em\hbox{[}}}
\def\cright{\hbox{]\kern-.16em\hbox{]}}}

\def\CC{ \mathbb{C} }
\def\DD{ \mathbb{D} }
\def\EE{ \mathbb{E} }

\def\into{ \rightarrowtail }
\def\onto{ \twoheadrightarrow }

\def\trio{ \triangleright}

\def\cleft{\hbox{[\kern-.16em\hbox{[}}}
\def\cright{\hbox{]\kern-.16em\hbox{]}}}

\theoremstyle{plain}
\newtheorem{dfn}[subsection]{Definition}
\newtheorem{thm}[subsection]{Theorem}
\newtheorem{lma}[subsection]{Lemma}
\newtheorem{prp}[subsection]{Proposition}
\newtheorem{cor}[subsection]{Corollary}

\theoremstyle{remark}

\newtheorem{exm}[subsection]{Example}

\begin{document}

\title{Partial Mal'tsevness and partial protomodularity}

\author{Dominique Bourn}

\begin{abstract}
We introduce the notion of Mal'tsev reflection which allows us to set up a partial notion of Mal'tsevness with respect to a class $\Sigma$ of split epimorphisms stable under pullback and containing the isomorphisms, and we investigate what is remaining of the properties of the global Mal'tsev context. We introduce also the notion of partial protomodularity in the non-pointed context.
\end{abstract}
 
\maketitle

\section*{Introduction}

A Mal'tsev category is a category in which any reflexive relation is an equivalence relation, see \cite{CLP} and \cite{CPP}.  The categories $Gp$ of groups and $K$-$Lie$ of Lie $K$-algebras are major examples of Mal'tsev categories. The terminology comes from the pioneering work of Mal'tsev in the varietal context \cite{Mat} which was later on widely developped in \cite{Sm}.

In \cite{B0}, Mal'tsev categories were characterized in terms of split epimorphisms: a finitely complete category $\DD$ is a Mal'tsev one if and only if any pullback of split epimorphisms in $\DD$:
$$
    \xymatrix{
    X'  \ar@<-4pt>[r]_(.6){\bar g} \ar@<-4pt>[d]_{f'} & X \ar@<-4pt>[l]_(.4){\bar t} \ar@<-4pt>[d]_f\\
    Y'   \ar@<-4pt>[r]_g \ar@<-4pt>[u]_{s'} & Y \ar@<-4pt>[l]_t \ar@<-4pt>[u]_s }
  $$
is such that the pair $(s',\bar t)$ is jointly extremally epic. 

More recently the same kind of property was observed, but only for a certain class $\Sigma$ of split epimorphisms $(f,s)$ which is stable under pullback and contains isomorphisms. By the classes of \emph{Schreier} or \emph{homogeneous split epimorphisms} in the categories $Mon$ of monoids and $SRng$ of semi-rings \cite{BMMS}, by the classes of \emph{puncturing} or \emph{acupuncturing split epimorphims} in the category of quandles \cite{B15}, by the class of \emph{split epimorphic functors with fibrant splittings} in the category $Cat_Y$ of categories with a fixed set of objects $Y$ \cite{B14}.

Starting from these observations, we are led to introduce a general concept of Mal'tsev reflection (Definition \ref{malref}) which allows us on the one hand to characterize the Mal'tsev categories as those categories $\mathbb C$ which are such that the fibred reflection of points $(\P_{\EE},I_{\EE})$ is a Mal'tsev reflection, and on the other hand to restrict the Mal'tsev concept to a subclass $\Sigma$ of split epimorphisms (Definition \ref{smal}). This amounts to say that the property about commutative squares of split epimorphisms described above is only demanded when the split epimorphism $(f,s)$ belongs to the class $\Sigma$ in question.

The aim of the articles is to investigate what does remain of the global Mal'tsev properties in this partial context. The Mal'tsev properties are classically dealing with three major themes: internal reflexive relations and categories, permutation of reflexive relations, centralization of reflexive relations. Roughly speaking one can say that, in any of these three domains, all the results remain valid provided that on the one hand we restrict our attention to reflexive relations or categories whose underlying reflexive graph:
$$\xymatrix{ X_1 \ar@<-6pt>[r]_{d_1} \ar@<6pt>[r]^{d_0} & X \ar[l]|{s_0} }$$
is \emph{such that the split epimorphism $(d_0,s_0)$ belongs to the class $\Sigma$} and on the other hand we restrict the basic slogan to: \emph{any such reflexive relation is transitive}, instead of: \emph{any reflexive relation is an equivalence relation}, see Sections \ref{1}, \ref{2} and \ref{3}. Another interesting point is that, since in the partial context it is no longer possible to rely upon the set theoretical argument of difunctionality of relations used in the original articles \cite{CLP} and \cite{CPP}, many proofs have had to be thoroughly renewed (see for instance Proposition \ref{permut}). 

The existence of preorders in $\Sigma$-Mal'tsev categories, see Examples \ref{order}, unexpectedly reveals that this weaker categorical structuration gives rise to an interesting flexibility in comparison with what can appear now as the rigidity of the global Mal'tsev structure.

When, in addition, the ground category is regular or efficiently regular we recover some aspects of the specific features of the regular Mal'tsev context as far as the construction of Baer sums of $\Sigma$-special extensions with a fixed abelian ``generalized kernel''. Finally we give a few words on the stronger notion of \emph{partial protomodularity} whose first aspects were investigated in \cite{BMMS} inside the stricter context of the pointed and point-congruous category $Mon$ of monoids.\\ 

\noindent The article is organized along the following lines:\\
\noindent Section 1 is devoted to the definition of the  \emph{partial Mal'tsev} context. Section 2 deals with the transitivity of the $\Sigma$-reflexive relations and the question of the permutation of reflexive relations. Section 3 is devoted to the question of the centralization of reflexive relations. Section 4 investigates the traces of another significant characterization of Mal'tsev categories  \cite{B0} given by base-change along split epimorphisms with respect to the fibration of points. Section 5 deals with the question of the existence of centralizers. Section 6 shows how a further left exact condition (point-congruousness) propagates the Mal'tsev properties and produces a Mal'tsev \emph{core}. Section 7 focuses on the regular context and describes the construction of the Baer sums of $\Sigma$-special extensions with abelian kernel relation. Section 8 introduces the notion of \emph{partial protomodularity} in the non-pointed context and precises what are the related consequences on the notion normal subobject; in this setting the Malt'sev core becomes a protomodular core.

\section{The fibration of points}

From now on, any category will be suppose finitely complete and split epimorphism will mean split epimorphism with a given splitting. Recall from \cite{B} that, for any category $\EE$, $\Pt(\EE)$ denotes the category whose objects are the split epimorphisms (=the ``genereralized points'') of $\EE$ and whose arrows are the commuting squares between such split epimorphisms, and that $\P_{\EE}:\Pt(\EE) \to \EE$ denotes the functor associating with each split epimorphism its codomain. 

This functor $\P_{\EE}:\Pt(\EE)\to\EE$ is a fibration (the so-called \emph{fibration of points}) whenever $\EE$ has pullbacks of split epimorphisms. The $\P_\EE$-cartesian maps are precisely pullbacks of split epimorphisms. Given any morphism $f:X\rightarrow Y$ in $\EE$, base-change along $f$ with respect to the fibration $\P_\EE$ is denoted by $f^*:\Pt_{Y}(\EE) \rightarrow \Pt_X(\EE)$. 

This fibration is underlying a reflection since it has a fully faithful right adjoint $I_{\EE}$ defined by $I_{\EE}(X)=(1_X,1_X)$. Given any reflection $(U,T):\mathbb C\rightleftarrows \mathbb D$, a map $\phi\in\mathbb C$ is said to be \emph{$U$-invertible} when $U(\phi)$ is invertible and \emph{$U$-cartesian} when the following diagram is a pullback in $\mathbb C$: 
$$\xymatrix@=12pt{
X  \ar[r]^{\phi} \ar[d]_{\eta_X} & Y   \ar[d]^{\eta_Y} \\
TU(X)  \ar[r]_{TU(\phi)}  & TU(Y)  } 
$$
This last definition is based upon the fact that this kind of maps is necessarily \emph{hypercartesian} with respect to the functor $U$.
\begin{dfn}\label{malref}
Let be given a reflection $(U,T):\mathbb C\rightleftarrows \mathbb D$. It will be said to be a Mal'tsev reflection when any square of split epimorphisms in $\mathbb C$:
$$\xymatrix{
X'  \ar@<-2pt>[r]_(.6){\bar g} \ar@<-2pt>[d]_{f'} & X  \ar@<-2pt>[l]_(.4){\bar t} \ar@<-2pt>[d]_{f} \\
Y'  \ar@<-2pt>[r]_{g} \ar@<-2pt>[u]_{s'} & Y \ar@<-2pt>[l]_{t} \ar@<-2pt>[u]_{s}} 
$$
where both $g$ and $\bar g$ are $U$-cartesian and both $f$ and $f'$ are $U$-invertible (which implies that this square is necessarily a pullback) is such that the pair $(s',\bar t)$ is jointly extremally epic.
\end{dfn}

\subsection{$\Sigma$-Mal'tsevness}

Let $\Sigma$ be a class of split epimorphisms. We shall denote by $\Sigma(\EE)$ the full subcategories of $Pt(\EE)$ whose objects are in the class $\Sigma$.

\begin{dfn}
The class $\Sigma$ is said to be:\\
1) fibrational when $\Sigma$ is stable under pullback and contains the isomorphisms\\
2) point-congruous when, in addition, $\Sigma(\EE)$ is stable under finite limits in $Pt(\EE)$.
\end{dfn}

The first point of the previous definition determines a reflexive subfibration of the fibration of points:
$$\xymatrix@=8pt{
 {\Sigma(\EE)\;} \ar@{>->}[rr]^j \ar[dd]_{\P^{\Sigma}_{\EE}} && Pt(\EE)  \ar[dd]_{\P_{\EE}} &&  \\
  &&&&&&\\
 \EE \ar@{=}[rr] \ar@<-1ex>@{.>}[uu]_{I^{\Sigma}_{\EE}} && \EE \ar@<-1ex>@{.>}[uu]_{I_{\EE}}
  }
  $$
where the reflection is defined   by $I^{\Sigma}_{\EE}(X)=(1_X,1_X)$, while the second point guarantees that any fibre $\Sigma_Y(\EE)$ is stable under finite limit in $Pt_Y(\EE)$ and any base-change $f^*:\Sigma_Y(\EE)\to \Sigma_X(\EE)$ is left exact.

\begin{dfn}\label{smal} 
Let $\Sigma$ be a fibrational class of split epimorphisms in $\EE$. Then $\EE$ is said to be a $\Sigma$-Mal'tsev category when the reflection $(\P^{\Sigma}_{\EE},I^{\Sigma}_{\EE})$ is a Mal'tsev one.
\end{dfn}
This amounts to say that the property about commutative squares of split epimorphisms described in the introduction is only demanded when the split epimorphism $(f,s)$ belongs to the class $\Sigma$. Accordingly \emph{a category $\CC$ is a Mal'tsev one if and only if the fibred reflection of points $(\P_{\EE},I_{\EE})$ is a Mal'tsev one}.

\subsection{Examples}
\textbf{1)} Let $Mon$ be the category of monoids. A split epimorphism $(f,s):X\rightleftarrows Y$ will be called a \emph{weakly Schreier} split epimorphism when, for any element $y\in Y$, the application $\mu_y:Ker f\to f^{-1}(y)$ defined by $\mu_y(k)=k\cdot s(y)$ is surjective. The class $\Sigma$ of weakly Schreier split epimorphisms is fibrational and  \emph{the category $Mon$ is a $\Sigma$-Mal'tsev category}. 
\proof Stability under pullback is straightforward. Given a pullback of split epimorphisms in $Mon$, as in the introduction, take any $(z,x)\in Y'\times_YX=X'$, i.e. such that $g(z)=f(x)$. Since $(f,s)$ is in $\Sigma$, there is some $k\in Ker f$ such that:
$$(z,x)=(z,k\cdot sf(x))=(z,k\cdot sg(z))=(1,k)\cdot(z,sg(z))=\bar t(k)\cdot s'(z)$$
So the only submonoid $W\subset Y'\times_YX=X'$ containing $s'(Y')$ and $\bar t(X)$ is $X'$.
\endproof
\noindent \textbf{1')}  In \cite{MMS} a split epimorphism $(f,s):X\rightleftarrows Y$ in $Mon$ was called a \emph{Schreier} split epimorphism when the application $\mu_y$ is bijective. This defines a sub-class $\Sigma'\subset \Sigma$ which was shown to be point-congruous in \cite{BMMS}; by Theorem 2.4.2 in this same article, the category $Mon$ is a $\Sigma'$-Mal'tsev category according to the present definition.

\smallskip
\noindent\textbf{2)}  Suppose that $U:\CC\to \DD$ is a left exact functor. It is clear that if $\Sigma$ is a fibrational (resp. point-congruous) class of split epimorphisms in $\DD$, so is the class $\bar{\Sigma}=U^{-1}{\Sigma}$ in $\CC$. When, in addition, the functor $U$ is conservative (i.e. reflects the isomorphisms), then $\CC$ is a $\bar{\Sigma}$-Mal'tsev category as soon as $\DD$ is a $\Sigma$-Mal'tsev one.

\smallskip
\noindent\textbf{3)}  Let $SRg$ be the category of semi-rings. The forgetful functor $U:SRg\to CoM$ (where $CoM$ is the category of commutative monoids) is left exact and conservative. We call \emph{weakly Schreier} a split epimorphism in $\bar{\Sigma}$. In \cite{BMMS} a split epimorphism in $\bar{\Sigma}'$ was called a Scheier one. Thanks to the point 2), this gives us two other partial Mal'tsev structures.

\smallskip
\noindent\textbf{4)}  A quandle is a set $X$ endowed with a binary idempotent operation $\trio :X\times X\to X$ such that for any object $x$ the translation $-\trio x:X\to X$ is an automorphism with respect to the binary operation $\trio$ whose inverse is denoted by $-\trio^{-1}x$. A homomorphism of quandles is an application $f:(X,\trio)\to (Y,\trio)$ which respects the binary operation. This defines the category $Qnd$ of quandles. The notion was independantly introduced in  \cite{jo} and \cite{Ma} in strong relationship with Knot Theory. In \cite{B15} a split epimorphism $(f,s):X\rightleftarrows Y$ in $Qnd$ was called a \emph{puncturing} (resp. \emph{acupuncturing}) split epimorphism when, for any element $y\in Y$, the application $s(y)\triangleleft-:f^{-1}(y)\to f^{-1}(y)$ is surjective (resp. bijective). The class $\Sigma$ of puncturing (resp. $\Sigma'$ of acupuncturing) split epimorphisms was shown to be fibrational (resp. point-congruous), and the category $Qnd$ was shown to be a $\Sigma$-Mal'tsev (and a fortiori a $\Sigma'$-Mal'tsev) category.

\smallskip
\noindent\textbf{5)}  The category $Mon$ is nothing but the fibre above the singleton $1$ of the fibration $U:Cat\to Set$ which associates with any category its set of objects. Let us denote by $Cat_Y$ the fibre above the set $Y$. In \cite{B14} a bijective on objects split epimorphic functor $(F,S):\mathbb Y\rightleftarrows \mathbb Y'$ in $Cat_Y$ was called \emph{with fibrant splittings} when any map $S(\phi)$ is cartesian. The class $\Sigma_Y$ of such split epimorphisms was shown to be point-congruous, and the category $Cat_Y$ to be a $\Sigma_Y$-Mal'tsev category.

\smallskip
\noindent\textbf{5')}  In the same article, a similar result was establish for the fibers of any fibration $(\;)_0:Cat\mathbb E\rightarrow \mathbb E$ where $Cat\mathbb E$ is the category of internal categories in the finitely complete category $\EE$.

\section{First Mal'tsev-type properties}

In this section we shall investigate the aspects of the global Mal'tsev properties which remain valid in the partial context of the $\Sigma$-Mal'tsevness. 

\subsection{$\Sigma$-relations and $\Sigma$-graphs}\label{1}

We recalled that in a Mal'tsev category any reflexive relation is an equivalence relation. Here we have to restrict our attention to the following:

\begin{dfn} \label{Sreflexiverel}(see also \cite{BMMS})
A graph $X_1$ (resp. a relation $R$) on an object $X$ will be said to be a $\Sigma$-graph ($\Sigma$-relation) when it is reflexive:
$$\xymatrix{ X_1 \ar@<-6pt>[r]_{d_1} \ar@<6pt>[r]^{d_0} & X \ar[l]|{s_0} }$$
and such that the split epimorphism $(d_0,s_0)$ belongs to the class $\Sigma$.
\end{dfn}
\noindent Our ground observation will be the following:
\begin{prp}\label{stranstit}
Let $\EE$ be a $\Sigma$-Mal'tsev category. Any (resp. symmetric) $\Sigma$-relation $S$ on an object $X$ is necessarily transitive (resp. an equivalence relation). 
\end{prp}
\proof
Let us recall that, given any reflexive relation $R$ on $X$ as on the right hand side below, its simplicial kernel is the upper part of the universal $2$-simplicial object associated with it:
$$
\xymatrix@=12pt
{
{K[d_0,d_1]\;} \ar@(u,u)[rr]^{\pi_{0}}\ar[rr]_{\pi_1} \ar@(d,d)[rr]_{\pi_{2}}  && R
\ar@<2ex>[rr]^{d_{0}} \ar@<-2ex>[rr]_{d_{1}} \ar@<3ex>[ll]^{\sigma_1}
\ar@<-2ex>[ll]_{\sigma_0}  && X \ar[ll]_{s_0} }
$$
When $\EE$ is finitely complete, $K[d_0,d_1]$ is obtained by the following pullback of reflexive graphs in $\EE$:
  $$\xymatrix@=12pt{
{K[d_0,d_1]\;} \ar[rrd]^{(\pi_0,\pi_1)}  \ar@<-1ex>[ddd]_{\pi_0} \ar@<1ex>[ddd]^{\pi_1}  \ar[rrrrr]^{\pi_2} &&&&&  {R\;} \ar[rrd]^{(d_0,d_1)}  \ar@<-1ex>[ddd]_{d_0}  \ar@<1ex>[ddd]^{d_1}\\
&& R[d_0] \ar@<-1ex>[ddll]_{d_0} \ar@<1ex>[ddll]^{d_1} \ar[rrrrr]^<<<<<<<<<{(d_1.d_0,d_1.d_1)}  &&&&& X\times X \ar@<-1ex>[ddll]_{p_0} \ar@<1ex>[ddll]^{p_1}  && \\
     &&& &&&&&&\\ 
R\ar[rrrrr]_{d_1}  \ar[uuu]_{} \ar[uurr]_{} &&&&& X \ar[uuu]_{} \ar[uurr]_{}
    }
    $$
   In $Set$-theoretical terms, $K[d_0,d_1]$ is the set of triple of elements $(x_0,x_1,x_2)\in X$ such that $x_0Rx_1Rx_2$ and $x_0Rx_2$.
The vertical part indexed by $0$ in the previous diagram  determines a factorization $(\pi_0,\pi_2):K[d_0,d_1]\rightarrow R\times_XR$ to the following vertical pullback:
$$
\xymatrix@=20pt{
{K[d_0,d_1]\;} \ar@<1ex>[rrd]^{\pi_2} \ar@<-1ex>[rdd]_{\pi_0} \ar@{>->}[rd]_{}\\
&  R\times_XR  \ar@<-4pt>[r]_(.6){d_2} \ar@<-4pt>[d]_{d_0} & R  \ar@<-2pt>[l]_(.4){s_1} \ar@<-4pt>[d]_{d_0} \ar@(u,u)@{.>}[ull]_{\sigma_1}\\
& R   \ar@<-4pt>[r]_{d_1} \ar@<-2pt>[u]_{s_0} \ar@(l,l)@{.>}[uul]^{\sigma_0} & X \ar@<-2pt>[l]_{s_0} \ar@<-2pt>[u]_{s_0} }
$$
In set theoretical terms, this factorization associates the pair $(x_0Rx_1,x_1Rx_2)$ with the triple $(x_0,x_1,x_2)$. It is a monomorphism since $(d_0,d_1): R\rightarrowtail X\times X$ is a relation. On the other hand, the dotted factorizations $\sigma_0$ and $\sigma_1$ complete the quadrangle into a commutative diagram of split epimorphisms. So, when $R$ is a $\Sigma$-relation, the factorization $(\pi_0,\pi_2):K[d_0,d_1]\rightarrow R\times_XR$ is an extremal epimorphism as well, and consequently an isomorphism. Accordingly the map $R\times_XR \stackrel{(\pi_0,\pi_2)^{-1}}{\longrightarrow} K[d_0,d_1] \stackrel{\pi_1}{\rightarrow} R$ produces the desired transitivity map.
\endproof

\begin{exm}\label{order} In a $\Sigma$-Mal'tsev category, there are reflexive relations which are not equivalence relations, and they are important ones:\\
1) The internal order in $Mon$ given by the usual order between natural
numbers:
\[ \xymatrix{ \mathcal{O}_{\mathbb{N}} \ar@<-6pt>[r]_{p_1} \ar@<6pt>[r]^{p_0} & \mathbb{N}, \ar[l]|{s_0} } \]
with: $\mathcal{O}_{\mathbb{N}} = \{ (x, y) \in \mathbb{N} \times \mathbb{N} \ | \ x \leq y\}$ is a Schreier-relation in $Mon$, see \cite{BMMS}; and a fortiori a weakly Schreier one.\\
2) Similarly the internal order in $Mon$ given by the usual order between integers on the group $\mathbb{Z}$:
\[ \xymatrix{ \mathcal{O}_{\mathbb{Z}} \ar@<-6pt>[r]_{p_1} \ar@<6pt>[r]^{p_0} & \mathbb{Z}, \ar[l]|{s_0} } \]
with: $\mathcal{O}_{\mathbb{Z}} = \{ (x, y) \in \mathbb{Z} \times \mathbb{Z} \ | \ x \leq y\}$ is a Schreier-relation in $Mon$, again see \cite{BMMS}.\\
3) More generally any partally ordered group $(G,\leq)$ supplies a similar example of Schreier-relation in $Mon$.
\end{exm}

In a Mal'tsev category, on a reflexive graph there is at most one structure of internal category, and it is actually an internal groupoid, see \cite{CPP} and \cite{Ja}. Here we get:

\begin{prp}\label{intcat}
Let $\EE$ be a $\Sigma$-Mal'tsev category. On a $\Sigma$-graph there is at most one structure of category. More precisely, it is sufficient to have a composition map: $d_1\colon X_2\rightarrow X_1$ where $X_2$ is the internal object of ``composable pair of morphisms'', which satisfies: $d_1s_0=1_{X_1}$ and $d_1s_1=1_{X_1}$.
\end{prp}
\proof
Starting with any $\Sigma$-graph: $\xymatrix{ X_1 \ar@<-6pt>[r]_{d_1} \ar@<6pt>[r]^{d_0} & X_0 \ar[l]|{s_0} }$,
let us consider the following pullback where $X_2$ (with simplicial indexations) is the internal object of ``composable pair of morphisms'':
$$
\xymatrix{
 X_2  \ar@<-4pt>[r]_(.6){d_2} \ar@<-4pt>[d]_{d_0} & X_1  \ar@<-4pt>[l]_(.4){s_1} \ar@<-4pt>[d]_{d_0} \\
 X_1   \ar@<-4pt>[r]_{d_1} \ar@<-4pt>[u]_{s_0} & X_0 \ar@<-4pt>[l]_{s_0} \ar@<-4pt>[u]_{s_0} }
$$
Since the right hand side split epimorphism is in $\Sigma$, the pair $(s_0,s_1):X_1\rightrightarrows X_2$ is jointly extremally epic. So, there is atmost one map $d_1:X_2\to X_1$ satisfying the Axioms 2 (composition with identities): $d_1s_0=1_{X_1}$, $d_1s_1=1_{X_1}$. The incidence Axioms 1: $d_0d_1=d_0d_0$, $d_1d_1=d_1d_2$ come for free by composition with the same jointly extremally epic pair.  In order to express the associativity we need the following pullback which defines $X_3$ as the internal objects of ``triples of composable morphims'' and where the split epimorphism $(d_2,s_1)$ is still in $\Sigma$:
$$\xymatrix{
 X_3  \ar@<-2pt>[r]_(.6){d_0} \ar@<-2pt>[d]_{d_3} & X_2  \ar@<-2pt>[l]_(.4){s_0} \ar@<-2pt>[d]_{d_2} \\
 X_2   \ar@<-2pt>[r]_{d_0} \ar@<-2pt>[u]_{s_2} & X_1 \ar@<-2pt>[l]_{s_0} \ar@<-2pt>[u]_{s_1} }
 $$
The composition map $d_1$ induces a unique couple of maps $(d_1,d_2)\colon X_3 \rightrightarrows X_2$ such that
$d_0d_1=d_0d_0$, $d_2d_1=d_1d_3$ and $d_0d_2=d_1d_0$, $d_2d_2=d_2d_3$. The associativity Axiom 3 is given by the remaining simplicial axiom: 
(3) $d_1d_1=d_1d_2$.
The checking of this axiom comes with composition with the pair $(s_0,s_2)$ of the previous diagram since it is jointly extremally epic as well.
\endproof
A reflexive graph with $d_0=d_1$ being just a split epimorphism, we get:
\begin{cor}\label{commut}
Let $\EE$ be a $\Sigma$-Mal'tsev category and $(f,s):X\rightleftarrows Y$ a split epimorphism in $\Sigma$. There is at most one structure of monoid on the object $(f,s)$ in the fibre $Pt_Y\EE$. When it is the case, this monoid is necessarily commutative.
\end{cor}
\proof
A monoid structure on the object $(f,s)$ in the fibre $Pt_Y\EE$ is just a category structure on the reflexive graph defined by $d_0=f=d_1$. Now consider the following pullback of split epimorphims:
$$
\xymatrix{
 R[f]  \ar@<-4pt>[r]_(.6){p_1} \ar@<-4pt>[d]_{p_0} & X  \ar@<-4pt>[l]_(.4){s_1} \ar@<-4pt>[d]_{f} \\
 X   \ar@<-4pt>[r]_{f} \ar@<-4pt>[u]_{s_{-1}} & Y \ar@<-4pt>[l]_{s} \ar@<-4pt>[u]_{s} }
$$
Let $m:R[f]\rightarrow X$ be the binary operation of the monoid and $tw:R[f]\rightarrow R[f]$ be the ``twisting isomorphism'' defined by $tw(x,x')=(x',x)$. Saying that the monoid is commutative is saying that $m.tw=m$ (*). Since $(f,s)$ is in $\Sigma$, the pair $(s_{-1},s_1)$ is jointly strongly epic and we check (*) by composition with this pair.
\endproof

\noindent We shall need further tools to get the characterization of internal groupoids, see section \ref{295}.

\subsection{Permutation of relations}\label{2}

One major aspect of Mal'tsev varieties \cite{Mat} and categories (provided they are regular \cite{CLP}) is that any pair of reflexive relations does permute. So, in this section, we shall suppose momentarily that $\EE$ is a \emph{regular category} \cite{Ba}. In any regular category, the binary relations can be composed. For that, consider two binary relations $R$ and $S$:
$$
    \xymatrix@=10pt{
    && T \ar@{.>}[dl]_{\pi_0^T} \ar@{.>}[dr]^{\pi_1^T} &&\\
     & R \ar[dl]_{d_0^R} \ar[dr]^<<{d_1^R} && S \ar[dl]_<<{d_0^S} \ar[dr]^{d_1^S} &\\
  X && Y && Z}
  $$
take the pullback of $d_0^S$ along $d_1^R$ and decompose the factorization $(d_0^R.p_0^T,d_1^S.p_1^T):T\to X\times Z$ into a monomorphism and a regular epimorphism: $T\onto R\circ S\rightarrowtail X\times Z$. Two relations $R$ and $S$ on an object $X$ are said to \emph{permute} when $R\circ S=S\circ R$. Recall that, if $R$ and $S$ is a pair of preorders (resp. equivalence relations) on $X$, then $R\circ S$ is a preorder (resp. an equivalence relation) when $R$ and $S$ permute. Moreover $R\circ S$ becomes the supremum of $R$ and $S$ among the preorders (resp. the equivalence relations). In a regular Mal'tsev category any pair of reflexive relations do permute \cite{CPP}, \cite{CLP}.
\begin{prp}\label{permut}
Let $\EE$ be a regular $\Sigma$-Mal'tsev category. Given any pair of a reflexive relation $R$ and a symmetric $\Sigma$-relation $S$ (and so an equivalence relation according to Proposition \ref{stranstit}) on a object $X$, the two relations permute.
\end{prp}
\proof
Since $S$ is a symmetric $\Sigma$-relation, both split epimorphisms $(d_0^S,s_0^S)$ and $(d_1^S,s_0^S)$ are in $\Sigma$. Let us denote by $R\square S$ the inverse image of the reflexive relation $S\times S$ along $(d_0^R,d_1^R):R\rightarrowtail Y\times Y$. This produces a double relation:
  $$ \xymatrix@=30pt{
      R \square S \ar@<-1,ex>[d]_{p_0^R}\ar@<+1,ex>[d]^{p_1^R} \ar@<-1,ex>[r]_{p_0^S}\ar@<+1,ex>[r]^{p_1^S}
     & S \ar@<-1,ex>[d]_{d_0^S}\ar@<+1,ex>[d]^{d_1^S} \ar[l]\\
       R \ar@<-1,ex>[r]_{d_0^R} \ar@<+1,ex>[r]^{d_1^R} \ar[u]_{} & X
  \ar[u]_{} \ar[l]
                     }
  $$
which is the largest double  relation on $X$ relating $R$ and $S$. In set theoretical terms, $R \square S$ defines the subset of elements $(u,v,u',v')$ of $X^4$ such that we have:
$$ \xymatrix@=10pt{
     u \ar@{.>}[r]^S \ar@{.>}[d]_R & v\ar@{.>}[d]^R\\
     u' \ar@{.>}[r]_S & v'
                     }
  $$
Now let us denote by $T$ the domain of the pullback of $d_0^S$ along $d_1^R$. Then the canonical factorization $\phi:R \square S\to T$ is a regular epimorphism, since the split epimorphism $(d_0^S,s_0^S)$ belongs to $\Sigma$. So, $R\circ S$ coincides with the decomposition of $(d_0^R.p_0^R,d_1^S.p_1^S):R \square S\to X\times X$ (which associates $(u,v')$ with the previous square) into a monomorphism and a regular epimorphism. If we denote by $\bar T$ the domain of the pullback of $d_0^R$ along $d_1^S$, the canonical factorization $\psi:R \square S\to\bar T$ is a regular epimorphism as well, since $(d_1^S,s_0^S)$ is in $\Sigma$. So, $S\circ R$ coincides with the decomposition of the morphism $(d_0^S.p_0^S,d_1^R.p_1^R):R \square S\to X\times X$ (which, again, associates $(u,v')$ with the previous square) into a monomorphism and a regular epimorphism. Since the maps $(d_0^R.p_0^R,d_1^S.p_1^S)$ and $(d_0^S.p_0^S,d_1^R.p_1^R)$ coincide, we get the permutation $R\circ S=S\circ R$.
\endproof 

\section{Commutation in the fibers $Pt_Y(\EE)$ and centralization of relations}\label{3}

The Mal'tsev setting was also shown to fit extremely well with the notion of centralization of equivalence relations, see \cite{Sm} for the varietal context and \cite{CPP}, \cite{Pe} and \cite{BG1} for the categorical one. In this section we shall investigate the question in the partial Mal'tsev setting. 

\subsection{Commutation in $Pt_Y(\EE)$}

Consider two maps having same codomain in the fibre $Pt_Y\EE$ and such that the split epimorphism $(f,s)$ is in $\Sigma$ as on the left hand side, and, as on the right hand side, consider the pullback of $f$ along $g$:
 $$
 \xymatrix@=20pt{
 && &&& & X'\ar[rd] \ar@{.>}[d]_{\phi} \ar[dl] &    \\
 {U\;} \ar[r]^{h}\ar@<-1ex>[dr]_{g} & V \ar[d]_{}  & {\;X} \ar[l]_{k} \ar[dl]_<<<<<{f} &&& {U\;} \ar[r]^{h}\ar@<-1ex>[dr]_{g} \ar@<1ex>[ru]^{s'} & V \ar[d]_{}  & {\;X} \ar[l]_{k} \ar[dl]_<<<<<{f} \ar@<-1ex>[ul]_{t'}\\
  & Y \ar[ul]_>>>>>>{t} \ar@<-1ex>[ru]_>>>>>{s} \ar@<-1ex>[u] & &&&   & Y \ar[ul]_>>>>>>{t} \ar@<-1ex>[ru]_>>>>>{s} \ar@<-1ex>[u]
              }
 $$
\begin{dfn}
 Let $\EE$ be a $\Sigma$-Mal'tsev category and $(f,s)$ a split epimorphism in $\Sigma$. The pair $(h,k)$ is said to commute in the fibre $Pt_Y\EE$ when there is a (necessarily unique) map $\phi:X'\rightarrow V$ such that $\phi.t'=k$ and $\phi.s'=h$. The map $\phi$ is called the \emph{cooperator} of this pair.
 \end{dfn}
 The unicity of $\phi$ comes from the fact that the pair $(t',s')$ in the right hand side diagram is jointly extremally epic and makes its existence a property for the pair $(h,k)$ and not a further data. Now we get the following significant characterization:
 \begin{prp}
 Let $\EE$ be a $\Sigma$-Mal'tsev category and
 \[ \xymatrix{ X_1 \ar@<-6pt>[r]_{d_1} \ar@<6pt>[r]^{d_0} & X_0 \ar[l]|{s_0} } \]
 a $\Sigma$-reflexive graph. It is an internal category if and only if the following subobjects:
 $$
 \xymatrix@=30pt{
 {X_1\;} \ar@{>->}[r]^-{(d_0,1_{X_{1}})}\ar@<-1ex>[rd]_{d_0} & X_0\times X_1 \ar@<-1ex>[d]_{p_{X_0}} & {\;X_1} \ar@{>->}[l]_-{( d_1,1_{X_{1}})} \ar@<2ex>[ld]^{d_1}\\
  & X_0 \ar[lu]_>>>>>>>{s_0} \ar[u]_-{(1_{X_{0}},s_0)} \ar@<-1ex>[ru]^>>>>>>{s_0}
  }
 $$
 do  commute in $Pt_{X_0}\EE$.
 \end{prp}
\proof
 The two subobjects commute in $Pt_{X_0}\EE$ if and only if they have a cooperator $\phi\colon X_2\rightarrow X_0\times X_1$, i.e\ a morphism satisfying $\phi. s_0=(d_1,1_{X_{1}})$ and $\phi. s_1=(d_0,1_{X_{1}})$:
 $$
 \xymatrix@=35pt{
  & X_2 \ar[d]_{\phi} \ar@<-1.5ex>[dl]_-{d_2} \ar@<1.5ex>[dr]^{d_0}\\
 {X_1} \ar@{>->}[r]^-{(d_0,1_{X_{1}})} \ar@<-2ex>[rd]_-{d_0} \ar@<.5ex>[ur]_>>>>>>{s_1} & X_0\times X_1 \ar@<-1ex>[d]_-{p_{X_0}} & {X_1} \ar@{>->}[l]_-{(d_1,1_{X_{1}})} \ar@<2ex>[ld]^-{d_1} \ar@<-.5ex>[ul]^>>>>>>{s_0}\\
  & X_0 \ar@<1ex>[lu]_>>>>>>>>{s_0} \ar[u]_{(1_{X_{0}},s_0)} \ar@<-1ex>[ru]^>>>>>>>{s_0}
  }
 $$
 where the whole quadrangle is the pullback which defines the internal object of composable pairs of the reflexive graph. So the morphism $\phi$ is necessarily a pair of the form $(d_0. d_2,d_1)$, where $d_1\colon X_2\rightarrow X_1$ is such that $d_1. s_0=1_{X_1},\; d_1. s_1=1_{X_1}$. Since the morphism $d_1$ satisfies these two identities, it makes the reflexive graph an internal category by Proposition \ref{intcat}. Conversely, the composition morphism $d_1\colon X_2\rightarrow X_1$ of an internal category satisfies the previous two identities and produces the cooperator $\phi=(d_0.d_2,d_1)$.
 \endproof
 
 \subsection{Centralization of reflexive relations}

Let $R$ and $S$ be two reflexive relations on an object $X$ in a category $\EE$. A \emph{double connecting relation} for this pair is a commutative diagram:
$$ \xymatrix@=30pt{
      W \ar@<-1,ex>[d]_{p_0^R}\ar@<+1,ex>[d]^{p_1^R} \ar@<-1,ex>[r]_{p_0^S}\ar@<+1,ex>[r]^{p_1^S}
     & S \ar@<-1,ex>[d]_{d_0^S}\ar@<+1,ex>[d]^{d_1^S} \ar[l]\\
       R \ar@<-1,ex>[r]_{d_0^R} \ar@<+1,ex>[r]^{d_1^R} \ar[u]_{} & X
  \ar[u]_{} \ar[l]
                     }
  $$
where the upper row and the left hand side vertical part are reflexive relations as well and the parallel pairs with same indexation $i\in \{0,1\}$ are underlying morphisms of reflexive relations. The diagram introducing $R\square S$ in the previous section produces the largest double connecting relation for the pair $(R,S)$.
\begin{dfn}
A centralization data on the pair $(R,S)$ is a double connecting relation such that the commutative square containing $d_0^S$ and $d_1^R$ is a pullback. We shall denote it in the following way:
$$ \xymatrix@=30pt{
      R\rtimes_XS \ar@<-1,ex>[d]_{p_0^R}\ar@<+1,ex>[d]^{p_1^R} \ar@<-1,ex>[r]_{p_0^S}\ar@<+1,ex>[r]^{p_1^S}
     & S \ar@<-1,ex>[d]_{d_0^S}\ar@<+1,ex>[d]^{d_1^S} \ar[l]\\
       R \ar@<-1,ex>[r]_{d_0^R} \ar@<+1,ex>[r]^{d_1^R} \ar[u]_{} & X
  \ar[u]_{} \ar[l]
                     }
  $$
When $S$ (resp. $R$) is a preorder, we demand in addition that the left hand side (resp. the upper) reflexive relation is so.
\end{dfn}
\noindent When $R$ and $S$ are equivalence relations, then this condition is equivalent to the fact that \emph{any of the commutative labelled squares is a pullback}. In set theoretical terms, to ask for a centralization data is equivalent (see \cite{BG1} and also \cite{K}, \cite{CPP}) to ask for an application $p:R\times_XS\to X$  satisfying:\\
1) the Mal'tsev conditions: $p(xRxSy)=y$ and $p(xRySy)=x$\\
2) the coherence conditions: $xSp(xRySz)$ and $p(xRySz)Rz$\\
3) when $S$ is a preorder, the left associativity condition:\\$(p(xRySz)RzSt)=p(xRySt)$\\
3')  when $R$ is a preorder, the right associativity condition:\\ $p(xRySp(yRzSt))=p(xRzSt)$.\\
Condition 2) means that, with any triple $xRySz$, we can associate a square of related elements:
  $$
  \xymatrix@=10pt {
    x \ar@{.>}[rr]^-{S} \ar[d]_-{R} & & {\;p(x,y,z)} \ar@{.>}[d]^-{R} \\
    y \ar[rr]_-{S} & & z. }
  $$
while conditions 1) give a coherence with the reflexivity and conditions 3) and 3') give a coherence with the transitivity.
The map $p$ underlying a centralization data is nothing but $p=d_1^S.p_0^S=d_0^R.p_1^R$. Recall that the notion of centralization data allows us to characterize the internal groupoids (again see \cite{K} and also \cite{CPP}, \cite{BG1}):
\begin{prp}\label{groupoid}
A reflexive graph in $\EE$: \[ \xymatrix{ X_1 \ar@<-6pt>[r]_{d_1} \ar@<6pt>[r]^{d_0} & X_0 \ar[l]|{s_0} } \] is underlying a groupoid structure if and only if there is a centralization data on the pair $(R[d_0],R[d_1])$.
\end{prp}
\proof
Any internal groupoid determines the following diagram  where, in set theoretical terms, $d_2$ is defined by $d_2(f,g)=g.f^{-1}$:
\begin{equation}\label{gr13}
  \xymatrix@=4pt{
{R[d_0]\;}   \ar@<-1ex>[ddd]_{d_0} \ar@<1ex>[ddd]^{d_1}  \ar[rrrrr]^{d_2} &&&&&  {X_1\;}  \ar@<-1ex>[ddd]_{d_0}  \ar@<1ex>[ddd]^{d_1}\\
&&&&\\
&&&&\\
X_1 \ar[rrrrr]_{d_1}  \ar[uuu]_{}  &&&&& X_0 \ar[uuu]_{} 
    } 
\end{equation}
and any of the commutative squares above is a pullback,. Complete this diagram by the kernel relations $R[d_1]$ and $R[d_2]$; this produces a centralization data for the pair $(R[d_0],R[d_1])$. 

We shall prove the converse, using the Yoneda lemma, by checking it in $Set$. The ternary map $p$ associated with a centralization data on the pair $(R[d_0],R[d_1])$ associates with any triple of arrows $(\alpha,\beta,\gamma)$ with the following incidence:
$$ \xymatrix@=20pt{
     u \ar@{.>}[rrd] \ar[d]_{\alpha} && v \ar[d]^{\gamma} \ar[dll]_<<<<{\beta} \\
     u'  && v'
                     }
  $$
a dotted arrow $\delta=p(\alpha,\beta,\gamma):u\to v'$. Starting with a composable pair $x\stackrel{f}{\rightarrow}y\stackrel{g}{\rightarrow}z$, set: $g.f=p(f,1_y,g)$.
\endproof
In this section we are going to show that, given a pair $(R,S)$ of a reflexive relation $(d_0^R,d_1^R): R\rightrightarrows X$ and a $\Sigma$-relation $S$ on an object $X$ in a $\Sigma$-Mal'tsev category $\EE$, there is at most one centralization data on this pair.

\begin{dfn}Let $\EE$ be a $\Sigma$-Mal'tsev category and $(R,S)$ a pair of a reflexive relation $R$ and a $\Sigma$-relation $S$ on the object $X$. We say that the two reflexive relations $R$ and $S$ centralize each other (which we shall denote by $[R,S]=0$ as usual) when the two following subobjects commute in the fibre $Pt_X(\EE)$:
$$
  \xymatrix@=10pt{
  {R\;} \ar@{>->}[rrr]^{(d_1^R,d_0^R)}\ar@<-1ex>[ddrrr]_{d_1^R} &&& X\times X \ar[dd]_{p_0}  &&& {\;S} \ar@{>->}[lll]_{(d_0^S,d_1^S)} \ar[ddlll]_{d_0^S} \\
  &&&&& \\
   && &X \ar[uulll]_{s_0^R} \ar@<-1ex>[rrruu]_{s_0^S} \ar@<-1ex>[uu]_{s_0} && &
               }
$$
\end{dfn}
Notice that we have $R^{op}$ on the left hand side. If we introduce the following pullback:
$$\xymatrix@=30pt{
      R \times_X S \ar@<-1,ex>[d]_{p_0^R} \ar@<+1,ex>[r]^{p_1^S}
     & S \ar@<-1,ex>[d]_{d_0^S} \ar[l]^{\sigma_0^S}\\
       R  \ar@<+1,ex>[r]^{d_1^R} \ar[u]_{\sigma_0^R} & X \ar[u]_{s_0^S} \ar[l]^{s_0^R}
                     }
$$
the cooperator is necessarily a pair $(\pi,p):R\times_XS\to X\times X$ where $\pi$ is $d_0^S.p_1^S=d_1^R.p_0^R$ and $p$ satisfies $p.\sigma_0^S=d_1^S$ and  $p.\sigma_0^R=d_0^R$, namely the Mal'tsev conditions (1). The map $p$ is called the \emph{connector} of the pair according to \cite{BG1}.

In a Mal'tsev category we know that two reflexive relations (i.e. two equivalence relations) $R$ and $S$ on $X$ centralize each other as soon as $R\cap S=\Delta_X$. Here in a $\Sigma$-Mal'tsev category, we have as well:
  \begin{prp}\label{cap}
  Let $\EE$ be a $\Sigma$-Mal'tsev category. The equivalence relation $R$ and the $\Sigma$-equivalence relation
  $S$ on the object $X$ centralize each other as soon as $R\cap S=\Delta_X$.
  \end{prp}
  \proof
 Consider the double relation:
  $$ \xymatrix@=30pt{
      R \square S \ar@<-1,ex>[d]_{p_0^R}\ar@<+1,ex>[d]^{p_1^R} \ar@<-1,ex>[r]_{p_0^S}\ar@<+1,ex>[r]^{p_1^S}
     & S \ar@<-1,ex>[d]_{d_0^S}\ar@<+1,ex>[d]^{d_1^S} \ar[l]\\
       R \ar@<-1,ex>[r]_{d_0^R} \ar@<+1,ex>[r]^{d_1^R} \ar[u]_{} & X
  \ar[u]_{} \ar[l]
                     }
  $$
  Since we have $R\cap S=\Delta_X$, the factorization $\theta: R\square S\rightarrow R\times_X S$ is necessarily a monomorphism. Since $S$ is a $\Sigma$-equivalence relation, the split epimorphism $(d_0^S,s_0^S)$ is in $\Sigma$, and this same factorization is an extremal epimorphism. Accordingly it is an isomorphism, which says that the square above which is downward indexed by $0$ and rightward indexed by $1$ is a pullback. Accordingly the map $R\square S \stackrel{p_0^S}{\rightarrow} S \stackrel{d_1^S}{\rightarrow} X$ produces the desired connector.
  \endproof
  
  \begin{prp}\label{double1}
  Let $\EE$ be a $\Sigma$-Mal'tsev category. Suppose the reflexive relation $R$ and the $\Sigma$-relation $S$ on $X$ centralize each other.  Then necessarily:\\
  1) the coherence conditions $xSp(xRySz)$ and $p(xRySz)Rz$ hold\\
  2) the left associativity $p(p(xRySz)RzSt)=p(xRySt)$ condition holds\\
  3) when $R$ is a preorder, the right associativity condition $p(xRySp(yRzSt))=p(xRzSt)$ holds as well.
  \end{prp}
  \proof
  Let us consider the following pullback:
   $$\xymatrix@=10pt {
      {U\;} \ar@{>->}[rr]^j \ar[d]_{} & & R\times_XS  \ar[d]^{(d_0^R.p_0^R,p)} \\
      {S\;} \ar@{>->}[rr]_{(d_0^S,d_1^S)} & & X\times X}
    $$
  It defines $U$ as the subobject of those $xRySz\in R\times_XS$ such that we have $xSp(xRySz)$. For any $ySz\in S$, the element $yRySz\in R\times_XS$ belongs to $U$ since we have $ySp(yRySz)$ by $p(yRySz)=z$. This means that $\sigma_0^S$ factors through $U$. In the same way, for any $xRy\in R$, the element $xRySy\in R\times_XS$ belongs to $U$ since we have $xSp(xRySy)$ by $p(xRySy)=x$. This means that $\sigma_0^R$ factors through $U$. Since the pair $(\sigma_0^R,\sigma_0^S)$ is jointly extremally epic, the map $j$ is an isomorphism, and for every $xRySz\in R\times_XS$ we have $xSp(xRySz)$.
  
  We have a similar result concerning the subobject $V\rightarrowtail R\times_XS$ defines by the following pullback:
  $$\xymatrix@=10pt {
    {V\;} \ar@{>->}[rr]^j \ar[d]_{} & & R\times_XS  \ar[d]^{(p,d_1^S.p_1^S)} \\
    {R\;} \ar@{>->}[rr]_{(d_0^R,d_1^R)} & & X\times X}
  $$
  This give us  $p(xRySz)Rz$ for any $xRySz\in R\times_XS$. So, both Mal'tsev conditions produce the following double connecting relation on the pair $(R,S)$ in $\EE$:
\begin{equation}
\xymatrix@=30pt{
      R \times_X S \ar@<-1,ex>[d]_{p_0^R}\ar@<+1,ex>[d]^{(p,d_1^S.p_1^S)} \ar@<-1,ex>[rr]_{(d_0^R.p_0^R,p)}\ar@<+1,ex>[rr]^{p_1^S}
     && S \ar@<-1,ex>[d]_{d_0^S}\ar@<+1,ex>[d]^{d_1^S} \ar[ll]\\
       R \ar@<-1,ex>[rr]_{d_0^R} \ar@<+1,ex>[rr]^{d_1^R} \ar[u]_{} && X
  \ar[u]_{} \ar[ll]
                     } 
\end{equation}
It is called the \emph{centralizing double relation} associated with the connector $p$ and it characterizes the centralization $[R,S]=0$. The left hand side vertical relation on $R$ can be describe in the following way: $(xRy)S_p(x'Ry')$ if and only if $ySy'$ and $p(xRySy')=x'$. 
Moreover, the square which is downward indexed by $0$ and rightward indexed by $1$ is a pullback. Accordingly since $S$ is a $\Sigma$-relation, so is the left hand side vertical relation which becomes a preorder. The pair $(d_0^R,(d_0^R.p_0^R,p))$ produces a morphism of preorder which means that $p(p(xRySz)RzSt)=p(xRySt)$.

Suppose, in addition, $R$ is a preorder as well. Consider the following diagram:
  $$ \xymatrix@=30pt{
  R^2\times_XS \ar@<-1,ex>[rr]_{p_0^S}\ar@{.>}[rr]_{}\ar@<+1,ex>[rr]^{p_2^S} \ar@<-1,ex>[d]_{p_0^{R^2}} &&  R \times_X S \ar@<-1,ex>[d]_{p_0^R} \ar@<-1,ex>[rr]_{(d_0^R.p_0^R,p)}\ar@<+1,ex>[rr]^{p_1^S}
     && S \ar@<-1,ex>[d]_{d_0^S} \ar[ll]\\
   R^2 \ar@<-1,ex>[rr]_{d_0^R} \ar@{.>}[rr]_{} \ar@<+1,ex>[rr]^{d_2^R} \ar[u]_{s_0^{R^2}}  &&  R \ar@<-1,ex>[rr]_{d_0^R} \ar@<+1,ex>[rr]^{d_1^R} \ar[u]_{s_0^R} && X \ar[u]_{s_0^S} \ar[ll]
                     }
  $$
where the dotted lower horizontal map is the transitivity morphism $d_1^R$ and the left hand side split epimorphism is the pullback of the middle one along $d_2^R$. This defines $R^2\times_XS$ as the elements $xRyRzSt$ and the map $p_2^S$ by $p_2^S(xRyRzSt)=yRzSt$. So, we have the factorization $p_0^S$ above $d_0^R$ defined by $p_0^S(xRyRzSt)=xRySp(yRzSt)$ and the dotted factorization $p_1^S$ above $d_1^R$ defined by $p_1^S(xRyRzSt)=xRzSt$. Saying that $p(xRySp(yRzSt))=p(xRzSt)$ is saying that $p_0^S$ and $p_1^S$ are coequalized by the map $(d_0^R.p_0^R,p)$. Now since $(d_0^S,s_0^S)$ is in $\Sigma$, the pair of sections induced by the pullback given by the whole rectangle is jointly extremally epic. The vertical section is $xRyRz\mapsto xRyRzSz$ and the horizontal one is $zSt\mapsto zRzRzSt$. The coequalization is question can be checked by composition with these two sections. Obviously we have: $p(xRySp(yRzSz))=p(xRySy)=x=p(xRzSz)$ and $p(zRzSp(zRzSt))=p(zRzSt)$.
\endproof
So, as expected, \emph{the previous double centralizing relation (2) produces an internal centralization data} for the pair $(R,S)$ and it is the unique possible one. The left hand side relation is an equivalence relation when so is $S$, and in this case the square which is downward indexed by $1$ and rightward indexed by $1$ is a pullback as well. The previous proposition shows that the upper row is an equivalence relation as soon as $R$ is so. In this case the square which is downward indexed by $0$ and rightward indexed by $0$ is a pullback as well. Finally when both $R$ and $S$ are equivalence relations, all the reflexive relations in this diagram are equivalence relations, and, moreover, any commutative square is a pullback.

We are now going to study the stability properties of the centralization. For that, we first need:

\begin{lma}\label{cartcart}
Let $\EE$ be a $\Sigma$-Mal'tsev category. Consider the following monomorphism between split epimorphisms in $Pt(\EE)$ where $(f,s)$ is in $\Sigma$:
$$\xymatrix@=4pt{
{\bar X\;}   \ar@<-1ex>[ddd]_{\bar f}  \ar@(u,u)@<1ex>[rrrrrrrrrr]^{\bar x} \ar@{>->}[rrrrr]^m &&&&&{X'\;}   \ar@<-1ex>[ddd]_{f'}  \ar@<1ex>[rrrrr]^{x} &&&&&  {X\;}  \ar@<-1ex>[ddd]_{f} \ar@(u,u)[llllllllll]^{\bar s_x} \ar[lllll]^{s_x} \\
&&&&&&&&\\
&&&&&&&&\\
{\bar Y\;} \ar@(d,d)@<1ex>[rrrrrrrrrr]^{\bar y}  \ar[uuu]_{\bar s} \ar@{>->}[rrrrr]_n &&&&& Y' \ar@<1ex>[rrrrr]^{y}  \ar[uuu]_{s'}  &&&&& Y \ar[uuu]_{s} \ar[lllll]^{s_y} \ar@(d,d)[llllllllll]^{\bar s_y}
    }
    $$
If the right hand side one is $\P_{\EE}$-cartesian (i.e. a pullback), so is the whole rectangle.
\end{lma}
\proof
Set $(\bar f',\bar s')=n^*(f,s)$, and denote $\theta:\bar X\to \bar X'$ the canonical factorization in $Pt_{\bar Y}\EE$. It is an extremal epimorphism since $(f,s)$ is in $\Sigma$. It is a monomorphism as well since $\bar n.\theta$ is the monomorphism $m$, where $\bar n=f'^*(n)$. Accordingly $\theta$ is an isomorphism.
\endproof
\begin{prp}\label{part1}
Let $\EE$ be a $\Sigma$-Mal'tsev category and $(R,S)$ be a pair of a reflexive relation and a $\Sigma$-relation $S$ on $X$ such that $[R,S]=0$. If $R'$ is a reflexive relation on $X$ such that $R'\subset R$, then we get $[R',S]=0$.
\end{prp}
\proof
Denote $i:R'\into R$ the inclusion, consider the reflexive relation $\bar R'=i^{-1}(R\rtimes_X S)$ on $R'$ and the following diagram:
$$\xymatrix@=4pt{
{\bar R'\;}   \ar@<-1ex>[ddd]_{p_0^{\bar R'}}  \ar@(u,u)@<1ex>[rrrrrrrrrr]^{\bar p_1^S} \ar@{>->}[rrrrr]^j &&&&&{R\rtimes_X S\;}   \ar@<-1ex>[ddd]_{p_0^R}  \ar@<1ex>[rrrrr]^{p_1^S} &&&&&  {S\;}  \ar@<-1ex>[ddd]_{d_0^S} \ar@(u,u)[llllllllll]^{\bar s_0^S} \ar[lllll]^{} \\
&&&&&&&&\\
&&&&&&&&\\
{R'\;} \ar@(d,d)@<1ex>[rrrrrrrrrr]^{p_1^{R'}}  \ar[uuu]_{} \ar@{>->}[rrrrr]_i &&&&& R \ar@<1ex>[rrrrr]^{d_1^R}  \ar[uuu]_{}  &&&&& X \ar[uuu]_{s_0^S} \ar[lllll]^{}\ar@(d,d)[llllllllll]^{s_0^{R'}}
    }
    $$
Then apply the previous lemma which produces the unique centralization data.
\endproof
\begin{prp}\label{invim}
Let $\EE$ be a $\Sigma$-Mal'tsev category. Let $(R,S)$ be a pair of reflexive relations on an object $X$ which is equipped with a centralization data $R\rtimes_XS$ and let $u:U\into X$ be a subobject such that $u^{-1}(S)$ is a $\Sigma$-relation. Then we have $[u^{-1}(R),u^{-1}(S)]=0$.
\end{prp}
\proof
Consider the map $d_1^S.p_0^S=d_0^R.p_1^R=p:R\rtimes_XS\to X$ and the following pullback of $p.\tilde{u}$ along $u$ where $\tilde u$ is the canonical factorization:
$$
  \xymatrix@=4pt{
{W\;}\ar@{>->}[rrrrr]^<<<<<<<<{\check u}\ar[dddrrrrr]^{\pi} &&&&&{u^{-1}(R)\times_Uu^{-1}(S)\;}    \ar@{>->}[rrrrr]^>>>>>>>>{\tilde u} &&&&&  R\rtimes_XS  \ar[ddd]^{p}\\
&&&&&&&&&\\
&&&&&&&&&&\\
&&&&& {U\;}   \ar@{>->}[rrrrr]_{u}  &&&&& X }
$$
In set theoretical terms $W=\{xRySz/ (x,y,z)\in U^3;\; p(xRySz)\in U\}$. Since $u^{-1}(S)$ is a $\Sigma$-relation, the following pair in jointly strongly epic:
$${u^{-1}(S)\;}\stackrel{s_0}{\into} u^{-1}(R)\times_Uu^{-1}(S)\stackrel{\sigma_0}{\leftarrowtail} {\;u^{-1}(R)}$$
The subobject $\check u$ contains $u^{-1}(S)$ since $p(xRxSy)=y \in U$ when $y$ is in $U$ and contains $u^{-1}(R)$ since $p(xRySy)=x \in U$ when $x$ is in $U$. Accordingly $\check u$ is an isomorphism, and $\pi.\check u^{-1}$ is the desired connector for the pair $(u^{-1}(R),u^{-1}(S))$.
\endproof

\subsection{Characterization of $\Sigma$-equivalence relations and internal $\Sigma$-groupoids}\label{295}

We recalled above that, when a reflexive graph:
\[ \xymatrix{ X_1 \ar@<-6pt>[r]_{d_1} \ar@<6pt>[r]^{d_0} & X_0 \ar[l]|{s_0} } \]
is underlying a groupoid, the diagram (1) is a discrete cofibration and a discrete fibration, namely that any of the commutative squares is a pullback.
 \begin{prp}\label{sigrel}
 Let $\EE$ be a $\Sigma$-Mal'tsev category. Consider a reflexive graph:
 \[ \xymatrix{ X_1 \ar@<-6pt>[r]_{d_1} \ar@<6pt>[r]^{d_0} & X_0 \ar[l]|{s_0} } \]
 The following conditions are equivalent:\\
 1) this graph is underlying a $\Sigma$-groupoid\\
 2) the kernel relation $R[d_0]$ is a $\Sigma$-relation and $[R[d_0],R[d_1]]=0$ \\
A reflexive relation $S$ on $X$ is a $\Sigma$-equivalence relation if and only if the map $d_0:S\to X$ is such that $R[d_0]$ is a $\Sigma$-relation, i.e. if and only if the morphism $d_0$ is $\Sigma$-special according to the terminology of Section \ref{sigspec}.
 \end{prp}
 \proof
If this graph is underlying a $\Sigma$-groupoid, this graph is a $\Sigma$-graph, and there is at most one such structure. Since the square indexed by $0$ in the diagram $(1)$ is a pullback, then $R[d_0]$ is a $\Sigma$-relation. Proposition \ref{groupoid} shows that $[R[d_0],R[d_1]]=0$.
Conversely suppose that $R[d_0]$ is a $\Sigma$-relation, then  $[R[d_0],R[d_1]]=0$ makes sense. This produces a double centralization data which, again by this same proposition, gives the groupoid structure.

Now suppose that the reflexive graph is actually a relation. The pair $(d_0,d_1):X_1\rightrightarrows X_0$ being jointly monic, we get $R[d_0]\cap R[d_1]=\Delta_{X_1}$. According to Proposition \ref{cap}, we have $[R[d_0],R[d_1]]=0$ as soon as $R[d_0]$ is a $\Sigma$-relation.
\endproof

\begin{cor}
Let $\EE$ be a $\Sigma$-Mal'tsev category and the following diagram a monomorphism of reflexif graphs:
$$
  \xymatrix@=4pt{
{U_1\;}   \ar@<-1ex>[ddd]_{d_0} \ar@<1ex>[ddd]^{d_1}  \ar@{>->}[rrrrr]^{u_1} &&&&&  {X_1\;}  \ar@<-1ex>[ddd]_{d_0}  \ar@<1ex>[ddd]^{d_1}\\
&&&&\\
&&&&\\
{U_0\;} \ar@{>->}[rrrrr]_{u_0}  \ar[uuu]_{}  &&&&& X_0 \ar[uuu]_{} 
    }
    $$
Suppose moreover that the left hand side vertical part is a $\Sigma$-graph. If the right hand side vertical diagram is underlying a groupoid structure, so is the left hand side vertical one.
\end{cor}
\proof
It is a straighforward application of Corollary \ref{invim}.
\endproof

\begin{cor}\label{abgp}
Let $\EE$ be a $\Sigma$-Mal'tsev category and $(f,s):X\rightleftarrows Y$ a split epimorphism. The following conditions are equivalent:\\
1) $R[f]$ is a $\Sigma$-relation and $[R[f],R[f]]=0$\\
2) the split epimorphism $(f,s)$ is in $\Sigma$ and is endowed with a (unique) group structure in $Pt_Y(\EE)$ which is necessarily abelian.\\
In this case, the split epimorphism is said to be \emph{abelian}.
\end{cor}
\proof
A split epimorphism is a reflexive graph with $d_0=d_1=f$. The commutativity of the group structure is guaranteed by Corollary \ref{commut}.
\endproof 
 
\section{Base-change along split epimorphisms}

In \cite{B0}, Mal'tsev categories were characterized by the fact that any base-change along a split epimorphism with respect to the fibration of points is fully faihtful and saturated on subobjects. We shall investigate in this section what does remain of this result in the partial context; this is given by Theorem \ref{saturated}. The whole section is based upon the following strong observation:
\begin{lma}\label{transal}
Let $\EE$ be a $\Sigma$-Mal'tsev category and the following quadrangle be any pullback of split epimorphisms:
$$
    \xymatrix{
    X'  \ar@<-4pt>[r]_(.6){\bar g} \ar@<-4pt>[d]_{f'} & X \ar@<-4pt>[l]_(.4){\bar t} \ar@<-4pt>[d]_f\\
    Y'   \ar@<-4pt>[r]_g \ar@<-4pt>[u]_{s'} & Y \ar@<-4pt>[l]_t \ar@<-4pt>[u]_s }
  $$
When its domain $(f',s')$ belongs to $\Sigma$, the upward and leftward square is a pushout.
\end{lma}
\proof
Consider any pair $(\phi,\sigma)$ of morphisms such that $\phi.s'=\sigma.g\;(*)$:
$$\xymatrix@=8pt{
        R[\bar g] \ar@<-1ex>[ddd]_{R(f')} \ar@<-1ex>[rrrr]_<<<<<<<<<{d_{1}^{x} }\ar@<1ex>[rrrr]^<<<<<<<<<{d_{0}^{x}} &&&&  {X'\;}  \ar[llll]  \ar@<-1ex>[ddd]_{f'}  \ar@{->>}[rrrr]^{\bar g}\ar[rrrrrd]_{\phi}  &&&&  {X\;}  \ar@<-1ex>[ddd]_{f}  \\
      &   &&&&   & & &  & T   && \\
      &&& &&&&&&\\
       R[g] \ar[uuu]_{R(s')} \ar@<-1ex>[rrrr]_{d_{1}^y} \ar@<1ex>[rrrr]^{d_{0}^y}  &&&&  Y' \ar[llll] \ar@{->>}[rrrr]_{g}  \ar[uuu]_{s'}
       && && Y   \ar@<-1ex>[uur]_{\sigma} \ar[uuu]_{s}
          }
          $$
and complete the diagram by the kernel relations $R[g]$ and $R[\bar g]$ which produce the left hand side pullbacks. The morphism $\bar g$ being a split epimorphism, it is the quotient of its kernel relation $R[\bar g]$. We shall obtain the desired factorization $X\to T$ by showing that $\phi$ coequalizes the pair $(d_0^x,d_1^x)$. The split epimorphism $(f',s')$ being in $\Sigma$, this can be done by composition with the jointly extremal pair $(R(s'),s_0^x); \; s_0^x:X'\into R[\bar g]$. This is trivial for the composition by $s_0^x$, and a consequence of the equality $(*)$ for the composition by $R(s')$.
\endproof

\subsection{Functors which are saturated on subobjects}

This is a technical section preparing the theorem concerning  the base-change functors in the next section. Recall that a subobject in a category $\EE$ is a class of monomorphisms up to isomorphism and a functor $F:\CC\to \DD$ is \emph{saturated on subobject} when it preserves monomophisms and produces a bijection between the set of subobjects of $X$ and  the set of subobjects of $F(X)$. We need now a refinement of this notion: suppose we have a commutative diagram of fully faithful subcategories:
$$\xymatrix@=10pt{
      {\CC'\;} \ar@{>->}[rrr] \ar[d]_{F'} &&& \CC \ar[d]^{F}\\
      {\DD'\;} \ar@{>->}[rrr] &&& \DD
        }
        $$
\begin{dfn}
The pair $(F,F')$ is said to be said to be fully faithful when, given any map $h:F(X)\to F(X')$ in $\DD$ with $X\in C'$ there is a unique $k:X\to X'$ in $\CC$ such that $F(k)=h$. The pair $(F,F')$ is said to be saturated on subobjects when $F$ preserves monomorphisms and produces a bijection  between the set of subobjects of $X$ with domain in $\CC'$ and  the set of subobjects of $F(X)$ with domain in $\DD'$.
\end{dfn}
\begin{lma}\label{sos}
Let the pair $(F,F')$ be saturated on subobjects and has a left inverse $(S,S')$. Then:\\
1) any monomorphism $m:A\rightarrowtail F(B)$ with $A\in \DD'$ is such that $m\simeq FS(m)$\\
2) any monomorphism $m:A\rightarrowtail F(B)$ with $A\in \DD'$ such that $S(m)$ is an isomorphism is itself an isomorphism\\
3) if, in addition, $\DD'$ is stable under finite limit and $F$ preserves products, any map $\phi:F(B)\rightarrow Z$ with $B\in \CC'$ and $Z\in \DD'$ such that $S(\phi)$ is a monomorphism is itself a monomorphism.
\end{lma}
\proof
1) Since $(F,F')$ is saturated on subobjects, starting with any monomorphism $m:A\rightarrowtail F(B)$ with $A\in \DD'$, there is a monomorphism $m':A' \rightarrowtail B$ with $A'\in \CC'$ such that $F(m')\simeq m$. Accordingly $m'=SF(m')\simeq S(m)$ and $m\simeq F(m')\simeq FS(m)$.\\
2) Accordingly, when $S(m)$ is an isomorphism, so is $m\simeq FS(m)$.\\
3) Let $\phi:F(B)\rightarrow Z$ be such that $B$ is in $\CC'$, $Z$ is in $\DD'$ and $S(\phi)$ is a monomorphism. Consider its kernel equivalence relation $(d_0^{\phi},d_1^{\phi}): R[\phi]\rightarrowtail F(B)\times F(B)=F(B\times B)$, the functor $F$ preserving products. The object $R[\phi]$ is in $\DD'$ since this subcategory is stable under finite limit. According to 1) we have $(d_0^{\phi},d_1^{\phi})\simeq FS(d_0^{\phi},d_1^{\phi})=(FS(d_0^{\phi}),FS(d_1^{\phi}))$. Since $S(\phi)$ is a monomorphism, we get $S(d_0^{\phi})=S(d_1^{\phi})$. Whence $d_0^{\phi}=d_1^{\phi}$; and so $\phi$ is a monomorphism.
\endproof

\subsection{Base-change along split epimorphisms}

We recalled that in a Mal'tsev category any base-change along a split epimorphism with respect to the fibration of points $\P_{\EE}$ is fully faithful and saturated on subobjects. We shall see here that, in the partial context, a similar property holds for the sub-fibration $\P^{\Sigma}_{\EE}$.

\begin{lma}\label{transv}
Let $\EE$ be a $\Sigma$-Mal'tsev category. Let be given any $\P_{\EE}$-cartesian split epimorphism as the quadrangle below. Then any other split epimorphism with its codomain $(f',s')$ in $\Sigma$:
$$\xymatrix@=15pt{
      {\bar X' \;} \ar@<-1ex>[ddd]_{\bar f'} \ar@{>->}[rd]_{\bar m} \ar@<1ex>[rrrr]^<<<<<<<<{\bar x} &&&&  {X'\;} \ar@{>->}[rd]^{m} \ar[llll]^>>>>>>>>{s_{\bar x}}  \ar@<-1ex>[ddd]_<<<<{f'} \\
    & \bar X \ar[ddl]_{\bar f}  \ar@<1ex>[rrrr]^<<<<{x} &&&&  X \ar[llll]^>>>>{s_x} \ar[ddl]_{f}   & & \\
    &&& &&&&&&\\
     \bar Y \ar[uuu]_{\bar s'}  \ar@<1ex>[rrrr]^{y} \ar@<-1ex>[uur]_{\bar s} &&&&  Y \ar[llll]^{s_y}   \ar[uuu]_>>>>{s'} \ar@<-1ex>[uur]_{s} 
        }
        $$
and with a monomorphic factorization $(m,\bar m)$ is necessarily $\P_{\EE}$-cartesian.
\end{lma}
\proof
Set $(\check f',\check s')=y^*(f',s')$ and denote by $\theta:\bar X'\to \check X'$ the induced factorization in the fibre $Pt_{\bar Y}\EE$. It is an extremal epimorphism since $(f',s')$ is in $\Sigma$. It is a monomorphism as well since we have $\theta.y^*(m)=\bar m$ which is a monomorphism. Accordingly $\theta$ is an isomorphism.
\endproof

\begin{thm}\label{saturated}
Suppose that $\EE$ is a $\Sigma$-Mal'tsev category and $(g,t):Y\rightleftarrows Z$ is any split epimorphism. Denote by $g_{\Sigma}^*:\Sigma_Z(\EE)\to \Sigma_Y(\EE)$ the restriction of the base-change. The pair $(g^*,g_{\Sigma}^*)$ is fully faithful and saturated on subobjects.
\end{thm}
\proof
  1) Full faithfulness. Consider the following diagram:
  $$\xymatrix@=10pt{
         {X'\;} \ar[rd]^{m'}   \ar@<-1ex>[ddd]_<<<<{f'}  \ar@{->>}[rrrr]^{g'}  &&&&  {X\;}  \ar@<-1ex>[ddd]_<<<<{f}  \\
      &   \bar X' \ar[ddl]_{\bar f'} \ar@{->>}[rrrr]^<<<<{\bar g}  & & &  & \bar X \ar[ddl]_{\bar f}  && \\
      &&& &&&&&&\\
     Y'  \ar@{->>}[rrrr]^{g}  \ar[uuu]_>>>>{s'} \ar@<-1ex>[uur]_{\bar s'} && && Y  \ar@<-1ex>[uur]_{\bar s} \ar[uuu]_>>>>{s}
          }
          $$
where the downward squares are pullback, $(f,s)$ is in $\Sigma$ and $m'$ a morphism in $Pt_{Y'}(\EE)$. Since $(f,s)$ is in $\Sigma$ and according to Lemma \ref{transal} the upward vertical square is a pushout; whence  a unique map $m:X\rightarrow \bar X$ such that $m.g'=\bar g.m'$ and $m.s=\bar s$; we get also $\bar f.m=f$ since $g'$ is a split epimorphism, and $m$ is a map in the fibre $Pt_Y(\EE)$ such that $g^*(m)=m'$.

\smallskip

2) Saturation on subobjects. First $g^*$ being left exact preserves monomorphisms. Consider now the following diagram where the right hand side quadrangle is a pullback, $(f',s')$ is in $\Sigma$ and $m$ is a monomorphism in $Pt_{Y}(\EE)$:
$$\xymatrix@=10pt{
      R[\bar g.m] \ar@<-1ex>[ddd]_{R(f')} \ar@{>->}[rd]_{R(m)} \ar@<-1ex>[rrrr]_<<<<<<<<{\delta_{1}}\ar@<1ex>[rrrr]^<<<<<<<<{\delta_{0}} &&&&  {X'\;} \ar@{>->}[rd]^{m} \ar[llll]  \ar@<-1ex>[ddd]_<<<<{f'}   &&&&   \\
    &  R[\bar g] \ar[ddl] \ar@<-1ex>[rrrr]_<<<<{d_{1}^{\bar g}} \ar@<1ex>[rrrr]^<<<<{d_{0}^{\bar g}} &&&&  \bar X' \ar[llll]
    \ar[ddl]_{\bar f'} \ar@{->>}[rrrr]^<<<<{\bar g}  & & &  & \bar X \ar[ddl]_{\bar f}  && \\
    &&& &&&&&&\\
     R[g] \ar[uuu]_{R(s')} \ar@<-1ex>[rrrr]_{d_{1}} \ar@<1ex>[rrrr]^{d_{0}} \ar@<-1ex>[uur] &&&&  Y' \ar[llll] \ar@{->>}[rrrr]^{g}  \ar[uuu]_>>>>{s'}
    \ar@<-1ex>[uur]_{\bar s'} && && Y   \ar@<-1ex>[uur]_{\bar s}
        }
        $$
Complete the diagram with the kernel relation $R[\bar g.m]$. According to Lemma \ref{transv}, any of the left hand side commutative squares is a pullback, and the following downward left hand side diagram is underlying a discrete fibration between equivalence relations:
$$\xymatrix@=10pt{
      R[\bar g.m] \ar@<-1ex>[ddd]_{R(f')}  \ar@<-1ex>[rrrr]_<<<<<<<<{\delta_{1}}\ar@<1ex>[rrrr]^<<<<<<<<{\delta_{0}} &&&&  {X'\;} \ar[llll] \ar@<3ex>[llll]^{\tau_1} \ar@<-1ex>[ddd]_{f'}  \ar@{.>>}[rrrr]^{g'} &&&&  {X\;}  \ar@<-1ex>[ddd]_{f} \ar@<1ex>[llll]^{\tau} \\
    &  &&&&   & & &  &  && \\
    &&& &&&&&&\\
     R[g] \ar[uuu]_{R(s')} \ar@<-1ex>[rrrr]_{d_{1}} \ar@<1ex>[rrrr]^{d_{0}}    &&&&  Y' \ar@<1ex>@(d,d)[llll]^{t_1} \ar[llll] \ar@{->>}[rrrr]^{g}  \ar[uuu]_{s'}
     && && Y \ar[uuu]_{s} \ar@<1ex>[llll]^{t}
        }
        $$
Let us denote by $t_1$ the unique map such that $d_1.t_1=1_Y$ and $d_0.t_1=t.g$. There is a unique map $\tau_1:X'\to R[\bar g.m]$ such that $\delta_1.\tau_1=1_{X'}$ and $R(f').\tau_1=t_1.f'$. It makes $(f',s')$ the pullback of $(R(f'),R(s'))$ along $t_1$. Now take the pullback $(f,s)$ of $(f',s')$ along $t$; it is in $\Sigma$ since so is $(f',s')$. Then it exists a unique $g':X'\to X$ which determines a pullback along $g$ and is such that $\tau.g'=\delta_0.\tau_1$. Moreover we can check that $g'.\delta_0=g'.\delta_1$ by composition with the extremally epic pair $(R(s'),\sigma_0)$. Accordingly the map $g'$ is the coequaliser of the pair $(\delta_0,\delta_1)$. In this way, the pair $(m,R(m))$ determines a factorization $n:X\to \bar X$ in the fibre $Pt_Y(\EE)$:
$$\xymatrix@=10pt{
      R[\bar g.m] \ar@<-1ex>[ddd]_{R(f')} \ar@{>->}[rd]_{R(m)} \ar@<-1ex>[rrrr]_<<<<<<<<{\delta_{1} }\ar@<1ex>[rrrr]^<<<<<<<<{\delta_{0}} &&&&  {X'\;} \ar@{>->}[rd]^{m} \ar[llll]  \ar@<-1ex>[ddd]_<<<<{f'}  \ar@{->>}[rrrr]^{g'} &&&&  {X\;}  \ar@<-1ex>[ddd]_{f}  \ar[rd]^{n}\\
    &  R[\bar g] \ar[ddl] \ar@<-1ex>[rrrr]_<<<<{d_{1}^{\bar g}} \ar@<1ex>[rrrr]^<<<<{d_{0}^{\bar g}} &&&&  \bar X' \ar[llll]
    \ar[ddl]_{\bar f'} \ar@{->>}[rrrr]^<<<<{\bar g}  & & &  & \bar X \ar[ddl]_{\bar f}  && \\
    &&& &&&&&&\\
     R[g] \ar[uuu]_{R(s')} \ar@<-1ex>[rrrr]_{d_{1}} \ar@<1ex>[rrrr]^{d_{0}} \ar@<-1ex>[uur] &&&&  Y' \ar[llll] \ar@{->>}[rrrr]_{g}  \ar[uuu]_>>>>{s'}
    \ar@<-1ex>[uur]_{\bar s'} && && Y \ar[uuu]_{s}  \ar@<-1ex>[uur]_{\bar s}
        }
        $$
  The upper right hand side quadrangle is a pullback since the two other right hand side commutative squares are so. Accordingly we get $m=g^*(n)$ and $n\simeq t^*g^*(m)$ is a monomorphism.
\endproof

\noindent We shall unwind now some consequences of the previous result which will be useful later on in the investigation about the relationship between $\Sigma$-Mal'tsevness and existence of centralizers, see Section \ref{existcentralizer}. For that we need the following:

\begin{dfn}(see also \cite{B0})
A reflexive graph in $Pt(\EE)$:
$$\xymatrix@=12pt
{
  X_1 \ar[dd]_{g_1} \ar@<2ex>[rr]^{d_{0}} \ar@<-2ex>[rr]_{d_{1}}& & X_0  \ar[dd]_-{g_0}  \ar[ll]_{s_0} \\
  &&& \\
  X_1' \ar@<2ex>[rr]^{d_0} \ar@<-2ex>[rr]_{d_1} \ar@<-1ex>[uu]_-{t_1}& & X_0'  \ar[ll]_{s_0} \ar@<-1ex>[uu]_-{t_0}
}
$$
will be said cartesian when any of the downward commutative squares is cartesian.
\end{dfn}

\begin{cor}
Let $\EE$ be a $\Sigma$-Mal'tsev category and the previous diagram be a cartesian split epimorphism of reflexive graphs. Then any subobject $(f_0,s_0)$ of $(g_0,t_0)$ in $Pt_{X'_0}(\EE)$ which is in $\Sigma$ is endowed with a unique up to isomorphism structure of subcartesian split epimorphism of reflexive graphs. When the upper graph of the diagram in question is underlying a category (resp. groupoid), the induced subgraph is a subcategory (resp. a subgroupoid). 
\end{cor}
\proof
Let $m:(f_0,s_0)\into (g_0,s_0)$ be the subobject in question. Now let the following upper reflexive graph be the fully faithful subgraph $m^{-1}(X_0,X_1)$ of the domain of the split epimorphism determined by the monomorphism $m$:
$$\xymatrix@=15pt{
      {M_1 \;} \ar@<-1ex>[ddd]_{f_1} \ar@{>->}[rd]_{m_1} \ar@<-1ex>[rrrr]_<<<<<<<<{d_{1}}\ar@<1ex>[rrrr]^<<<<<<<<{d_{0}} &&&&  {M_0\;} \ar@{>->}[rd]^{m} \ar[llll]  \ar@<-1ex>[ddd]_<<<<{f_0} \\
    &  X_1 \ar[ddl] \ar@<-1ex>[rrrr]_<<<<{d_{1}} \ar@<1ex>[rrrr]^<<<<{d_{0}} &&&&  X_0 \ar[llll]
    \ar[ddl]_{g_0}   & & \\
    &&& &&&&&&\\
     X'_1 \ar@{.>}[uuu]_{s_1} \ar@<-1ex>[rrrr]_{d_{1}} \ar@<1ex>[rrrr]^{d_{0}} \ar@<-1ex>[uur]_{t_1} &&&&  X'_0 \ar[llll]   \ar[uuu]_>>>>{s_0}
    \ar@<-1ex>[uur]_{t_0} 
        }
        $$
So there is a factorization $s_1$ which completes the vertical part as a square of split epimorphisms. According to the previous lemma and since $(f_0,s_0)$ is in $\Sigma$, this vertical part is made of pullbacks and consequently gives rise to a cartesian split epimorphism of reflexive graphs. As a full subgraph of $(X_0,X_1)$, the graph $(M_0,M_1)$ is a subcategory (resp. a subgroupoid) as soon as $(X_0,X_1)$ is a category (resp. a groupoid). 
\endproof

\begin{cor}
Let $\EE$ be a $\Sigma$-Mal'tsev category and $(R,T)$ a pair of reflexive relations on $X$ equipped with a centralization data. Suppose $S$ is a $\Sigma$-relation such that $S\subset T$. Then we have $[R,S]=0$.
\end{cor}
\proof
The split epimorphism $(d_0^S,s_0^S)$ of $(d_0^T,s_0^T)$ is in $\Sigma$. According to the previous corollary, the centralization data on $(R,T)$ produces a cartesian split epimorphism of reflexive graphs:
$$\xymatrix@=15pt
{
  R_1 \ar[dd]_-{p_0^R} \ar@<1ex>[rr]^{p_{0}^S} \ar@<-1ex>[rr]_{p_{1}^S}& & S  \ar[dd]_-{d_0^S}  \ar[ll]_{} \ar@(r,r)@{.>}[dd]^{d_1^S} \\
  &&& \\
  R \ar@<1ex>[rr]^{d_0^R} \ar@<-1ex>[rr]_{d_1^R} \ar@<-1ex>[uu]_-{}& & X  \ar[ll]_{} \ar@<-1ex>[uu]_-{s_0^S}
}
$$
The map $p=d_1^S.p_0^S$ produces the desired connector for the pair $(R,S)$.
\endproof

So we finally got here the second part of what was a symmetric result in the global Mal'tsev context and of which the first part is given by Proposition \ref{part1}. 

\section{Existence of centralizers}\label{existcentralizer}

\noindent In this section, we shall characterize the existence of centralizers in the $\Sigma$-Mal'tsev categories in the same way as in the Mal'tsev context, see \cite{B13}. For that, we shall be interested in the cartesian reflexive relations in $\Sigma(\EE)$ and shall begin with the following observation:
 
 \begin{lma}\label{lem1}
 Let $\EE$ be a $\Sigma$-Mal'tsev category and $(f,s)$ a split epimorphism in $\Sigma$. Let $R$ be a cartesian reflexive relation on  $(f,s)$ in $Pt(\EE)$. Then any reflexive sub-relation $j:T\rightarrowtail R$ of $R$ is cartesian as well.
 \end{lma}
 \proof
 Consider the following diagram:
 $$ \xymatrix@=18pt{
   {T_X\;} \ar[d]_{T_{f}}\ar@<1ex>@(u,u)[rrrr]^{d_i^X} \ar@{>->}[rr]^{j_X}&&{R_X\;} \ar[d]_{R_f}   \ar@<1ex>[rr]^{d_i^X} && {\;X} \ar[d]_{f} \ar[ll] \ar@<1ex>@(u,u)[llll]\\
   {T_Y\;} \ar@{>->}[rr]_{j_Y}\ar@<-1ex>@(d,d)[rrrr]_{d_i^Y} \ar@<-1ex>[u]_{T_{s}} && {R_Y\;}  \ar@<-1ex>[u]_{R_s}  \ar@<1ex>[rr]^{d_i^Y} & & {\;Y}  \ar@<-1ex>[u]_{s} \ar@{>->}[ll] \ar@(d,d)[llll]
                     }
  $$
And apply Lemma \ref{cartcart} for $i\in \{0,1\}$.
\endproof
 The cartesian reflexive relations allow us to characterize the centralization of a pair $(R,S)$:
 \begin{lma}\label{lem2}
 Let $\EE$ be a $\Sigma$-Mal'tsev category and $S$ a $\Sigma$-relation on an object $X$. Given any equivalence relation $R$ on $X$, we have $[R,S]=0$ if and only $R$ is underlying a cartesian equivalence relation on the split epimorphism $(d_0^S,s_0^S)$.
 \end{lma}
 \proof
 If we have $[R,S]=0$,consider the associated double centralizing relation: 
 $$ \xymatrix@=35pt{
       R \times_X S \ar@<-1,ex>[d]_{p_0^R}\ar@<+1,ex>@{.>}[d] \ar@<-1,ex>[r]_{(d_0^R.p_0^R,p)}\ar@<+1,ex>[r]^{p_1^S}
      & S \ar@<-1,ex>[d]_{d_0^S}\ar@<+1,ex>@{.>}[d]^{d_1^S} \ar[l]\\
        R \ar@<-1,ex>[r]_{d_0^R} \ar@<+1,ex>[r]^{d_1^R} \ar[u]_{} & X
   \ar[u]_{} \ar[l]
                      }
   $$
Its non-dotted part provides us with a cartesian equivalence relation above $R$ on the split epimorphism $(d_0^S,s_0^S)$. Conversely suppose that we have a cartesian reflexive relation above $R$ on $(d_0^S,s_0^S)$:
   $$ \xymatrix@=30pt{
       \bar R \ar@<-1ex>[d]_{\delta_0^R}  \ar@<-1ex>[rr]_{\delta_0^S}\ar@<1ex>[rr]^{\delta_1^S} && S \ar@<-1ex>[d]_{d_0^S}\ar@(r,r)@{.>}[d]^{d_1^S}  \ar[ll]\\
       R \ar@<-1ex>[rr]_{d_0^R} \ar@<1ex>[rr]^{d_1^R} \ar[u]_{} && Y  \ar[u]_{s_0^S} \ar[ll] 
                      }
   $$
 The rightward square with the $d_1$ being a pullback, the object $\bar R$ is $R\times_XS$ and the map $p=d_1^S.\delta_0^S$ provides us the desired connector. \endproof
 In a $\Sigma$-Mal'tsev category $\EE$, we can introduce, in the same way as in a Mal'tsev context, the following:
 
 \begin{dfn}
 Let $\EE$ be a $\Sigma$-Mal'tsev category. The centralizer of a $\Sigma$-equivalence relation $S$ on an object $X$ is the largest equivalence relation on $X$ centralizing $S$.
 \end{dfn}
 Similarly to the Mal'tsev setting \cite{B13}, the existence of centralizers of $\Sigma$-equivalence relations can be characterized by a property of the fibration $\P_{\EE}^{\Sigma}$. For that let us introduce the following:
 \begin{dfn}
 We say that a $\Sigma$-Mal'tsev category $\EE$ is (resp. strongly) action distinctive when, on any split  $\Sigma$-special epimorphism $(f,s)$ (resp.  any split epimorphism $(f,s)$ in $\Sigma$), there exists a largest cartesian equivalence relation on it which we shall denote in the following way:
 $$ \xymatrix@=20pt{
     D_X[f,s] \ar@<-1,ex>[d]_{D_f} \ar@<-1ex>[r]_{\delta_1^X}\ar@<1ex>[r]^{\delta_0^X} & X \ar[d]_{f} \ar[l]\\
     D_Y[f,s]  \ar@<-1ex>[r]_{\delta_1^Y} \ar@<1ex>[r]^{\delta_0^Y} \ar[u]_{D_s} & Y  \ar@<-1ex>[u]_{s} \ar[l] 
                    }
 $$
 and shall call the action distinctive equivalence relation on $(f,s)$.
 \end{dfn}
 In these contexts, according to Lemma \ref{lem1}, an equivalence relation $R$ on a split $\Sigma$-special epimorphism (resp. on a split epimorphism $(f,s)$ in $\Sigma$) is cartesian if and only if it is a smaller than $\mathbb D[f,s]$.
  
 \begin{thm}\label{zent}
 Let $\EE$ be a $\Sigma$-Mal'tsev category. It is action distinctive if and only if any $\Sigma$-equivalence relation $(d_0^S,d_1^S):S\rightrightarrows X$ admits a centralizer.
 \end{thm}
 \proof 
 Suppose the $\Sigma$-Mal'tsev category $\EE$ is action distinctive. According to Proposition \ref{sigrel}, an equivalence relation $S$ is a $\Sigma$-equivalence relation if and only if $d_0^S$ is a $\Sigma$-special split epimorphism. So we can consider the action distinctive equivalence relation on  $(d_0^S,s_0^S)$
 $$ \xymatrix@=20pt{
     D_R[d_0,s_0] \ar@<-1ex>[d]_{D_{d_0}} \ar@<-1ex>[r]_{\delta_1^R}\ar@<1ex>[r]^{\delta_0^R} & S \ar@<-1ex>[d]_{d_0^S} \ar[l]\\
     D_X[d_0,s_0] \ar@<-1ex>[r]_{\delta_1^X} \ar@<1ex>[r]^{\delta_0^X} \ar[u]_{D_{s_0}} & X  \ar[u]_{s_0^S} \ar[l] 
                    }
 $$
 According to Lemma \ref{lem2}, we have $[D_X[d_0,s_0],R]=0$. Let us show that the equivalence relation $D_X[d_0,s_0]$ is the centralizer of $S$. Let $R$ be any equivalence relation on $X$ such that $[R,S]=0$. According to this same lemma, it determines a cartesian equivalence relation on $(d_0^S,s_0^S)$ above $R$, whence $R\subset D_X[d_0,s_0]$ since $D_X[d_0,s_0]$ is underlying the largest cartesian one.
 
 Conversely, suppose the $\Sigma$-Mal'tsev category $\EE$ has centralizers of $\Sigma$-equiva-\linebreak lence relations. Let $(f,s)$ be a $\Sigma$-special split epimorphism; then its kernel relation $R[f]$ is a $\Sigma$-equivalence relation which admits a centralizer $Z[R[f]]$:
 $$\xymatrix@=20pt
 {
 \Sigma_f \ar@<-1ex>[d]_{\delta_0^Z} \ar@<1ex>@{.>}[d]^{\delta_1^Z} \ar@<-1ex>[rr]_{\delta_1^R}\ar@<1ex>[rr]^{\delta_0^R} && R[f] \ar@<-1ex>[d]_{d_0^f} \ar@<1ex>@{.>}[d]^{d_1^f} \ar[ll]\\
 Z[R[f]] \ar@<-1ex>[rr]_{d_1^Z} \ar@<1ex>[rr]^{d_0^Z} \ar[u]_{} && X  \ar[u] \ar[ll] 
 }
 $$
 The non-dotted part of the previous diagram determines a cartesian equivalence relation on the split epimorphism $(d_0^f,s_0^f)$, and according to the characterization given by Lemma \ref{lem2}, it is certainly the largest one. Accordingly, it is the action distinctive equivalence relation on $(d_0^f,s_0^f)$ and we shall denote it by $\mathbb D[d_0^f,s_0^f]$.
 Now consider the following  $\P_{\EE}^{\Sigma}$-cartesian monomorphism:
 $$\xymatrix@=15pt
 {
 X   \ar@<-1ex>[d]_{f}\ar[r]^{s_1} & R[f] \ar@<-1ex>[d]_{d_{0}^f} \\
 Y \ar[r]_{s} \ar[u]_{s}  & X \ar[u]_{s_0^f}  
 }
 $$
 Define the desired equivalence relation $\mathbb D[f,s]$ as $(s,s_1)^{-1}(\mathbb D[d_0^f,s_0^f])$. It is clearly a cartesian equivalence relation as an inverse image of cartesian equivalence relation along a cartesian map. Let us show it has the required universal property. So, consider any cartesian equivalence relation $(R,T)$ on $(f,s)$, and complete the following lower pullbacks with the kernel equivalence relations:
 $$ \xymatrix@=20pt{
     R[g] \ar@<-1,ex>[d]_{d_0}\ar@<+1,ex>[d]^{d_1} \ar@<-1,ex>[r]_{R(d_0)}\ar@<+1,ex>[r]^{R(d_1)} & R[f] \ar[l] \ar@<-1,ex>[d]_{d_0^f}\ar@<+1,ex>[d]^{d_1^f} \\
      T \ar@<-1,ex>[r]_{d_{0}^T} \ar@<+1,ex>[r]^{d_{1}^T} \ar[u]_{} \ar@<-1ex>[d]_{g} & X \ar[l] \ar@<-1ex>[d]_{f} \ar[u] \\
      R \ar@<-1,ex>[r]_{d_{0}^R} \ar@<+1,ex>[r]^{d_{1}^R} \ar[u]_{t}  & Y  \ar[l] \ar[u]_{s}
                    }
 $$
 First, notice that $R=s^{-1}(T)$. Since the lower diagrams are pullbacks, so are the upper ones which determine a double centralizing relation implying that we have $[T,R[f]]=0$. Accordingly, we get an inclusion $j: T\rightarrowtail Z[R[f]]$ and consequently $R=s^{-1}(T)\subset s^{-1}(Z[R[f]])=D_Y[f,s]$ by the definition of $\mathbb D[f,s]$ above. Finally consider the following diagram where the map $\check j$ is the unique factorization induced by the map $d_0^T$:
 $$
 \xymatrix@=20pt{
  {T\;} \ar@{>->}[r]^<<<<<<{\check j} \ar@<-1ex>[d]_{g} \ar@(u,u)[rr]^{d_0^T} & {s_1^{-1}(\Sigma_f)\;} \ar@<-1,ex>[r]_{d_{0}} \ar@<+1,ex>@{.>}[r]^{d_{1}} \ar@<-1ex>[d]  & {\;X} \ar[d]_{f} \ar@{>->}[l] \\
 {S\;} \ar@{>->}[r]^<<<<{j} \ar[u]_{t} \ar@(d,d)[rr]_{d_0^S}  & {s^{-1}(Z[R[f]])\;} \ar@<-1,ex>[r]_{d_{0}} \ar@<+1,ex>@{.>}[r]^{d_{1}}\ar[u]_{}  & {\;Y} \ar@<-1ex>[u]_{s} \ar@{>->}[l]
 }             
 $$
 To get $T\subset s_1^{-1}(\Sigma_f)=D_X[f,s]$, it remains to show that we have $d_1.\check j=d_1^T$; but, the split epimorphism $(f,s)$ being in $\Sigma$ as a split $\Sigma$-special epimorphism, the pair $(t,s_0^T)$, with $s_0^T:X\rightarrow T$, is jointly strongly epic, and this equality can be checked by composition with this pair.
\endproof

\begin{exm}
\textbf{1)} The category $Mon$ of monoids is an action distinctive $\Sigma'$-Mal'tsev category, Section 5.4 in \cite{BMMS}.\\
\textbf{2)} The category $SRg$ of semi-rings is an action distinctive $\bar{\Sigma}'$-Mal'tsev category, Section 6.6 in \cite{BMMS}.\\
\textbf{3)} Any fibre $Cat_Y$ is an action distinctive $\Sigma_Y$-Mal'tsev category, Section 6.1 in \cite{B14}.
\end{exm}

\section{Point-congruous $\Sigma$-Mal'tsev categories}

In this section we shall show that, when in addition the class $\Sigma$ is point-congruous, there are very important consequences for the category $\EE$; this produces, in particular, a Mal'tsev \emph{core}.

\subsection{Reflexive graphs in $\Sigma(\EE)$}

In a Mal'tsev category we have also the following result, see \cite{B0}: given any split epimorphism of reflexive graphs,
$$\xymatrix@=12pt
{
  X_1 \ar[dd]_-{g_1} \ar@<2ex>[rr]^{d_{0}} \ar@<-2ex>[rr]_{d_{1}}& & X_0  \ar[dd]_-{g_0}  \ar[ll]_{s_0} \\
  &&& \\
  X_1' \ar@<2ex>[rr]^{d_0'} \ar@<-2ex>[rr]_{d_1'} \ar@<-1ex>[uu]_-{t_1}& & X_0'  \ar[ll]_{s_0'} \ar@<-1ex>[uu]_-{t_0}
}
$$
the commutative square with maps $d_1$ is a pullback as soon as so is the square with maps $d_0$; in other words this split epimorphism of reflexive graphs is cartesian as soon as the square with maps $d_0$ is a pullback. Here, similarly, we have:
\begin{prp}\label{splitfib}
Let $\EE$ be a point-congruous $\Sigma$-Mal'tsev category. Let be given a split epimorphism of reflexive graphs, with the split epimorphism $(g_0,t_0)$ in $\Sigma$. It is cartesian as soon as the square with maps $d_0$ is a pullback.
\end{prp}
\proof
Set $(\bar g,\bar t)=d_1'^*(g_0,t_0)$; it is in $\Sigma$ since so is $(g_0,t_0)$. Denote $\theta:X_1\to \bar X$ the induced factorization in $Pt_{X_1'}(\EE)$. It is an extremal epimorphism since $(g_0,t_0)$ is in $\Sigma$, and actually it lies in $\Sigma_{X'}(\EE)$ with the assumption about $d_0$. Now the leftward square is a pullback, and to say it is equivalent to say that $s_0^*(\theta)=1_X$.  By Theorem \ref{saturated} applied to the base-change $d_0^{'*}$ and the point 3) in Lemma \ref{sos} which holds when $\Sigma$ is point-congruous, this morphism $\theta$ is a monomorphism since $s_0^*(\theta)=1_X$. Accordingly it is an isomorphism.
\endproof

\subsection{$\Sigma$-special morphisms}\label{sigspec}

\begin{dfn}
Let $\EE$ be a $\Sigma$-Mal'tsev category. A morphism $f:X\to Y$ is said to be $\Sigma$-special when its kernel relation $R[f]$ is a $\Sigma$-relation. An objet $X$ is said to be $\Sigma$-special when its terminal map $\tau_X:X\to 1$ is $\Sigma$-special.
\end{dfn}

When the class $\Sigma$ is fibrational, i.e. stable under pullback, so is the class of $\Sigma$-special morphisms.  It is also clear that any isomorphism is $\Sigma$-special. A split epimorphism which is in $\Sigma$ is not necessarily $\Sigma$-special; however we have the following:
\begin{lma}\label{splitsp}
Let $\EE$ be $\Sigma$-Mal'tsev category. Any $\Sigma$-special split epimorphism $(f,s)$ is in $\Sigma$.
\end{lma}
\proof
It is a consequence of the fact that the following leftward square is a pullback:
$$\xymatrix@=12pt
{
  R[f] \ar[dd]_{p_0}  & & X  \ar[dd]_{f}  \ar[ll]_{s_1} \\
  &&& \\
  X  \ar@<-1ex>[uu]_{s_0} & & Y  \ar[ll]^{s} \ar@<-1ex>[uu]_{s}
}
$$
for any split epimorphism $(f,s)$.
\endproof

We shall denote by $\Sigma l(\EE)$ the category whose objects are the $\Sigma$-special morphisms and whose morphisms are the commutative squares between them. When $\Sigma$ is point-congruous, $\Sigma l(\EE)$ is stable under finite limit in $\EE^2$. Similarly we shall denote by $\Sigma l_Y\EE$ the full subcategory of the slice category $\EE/Y$ whose objects are the $\Sigma$-special morphisms.
\begin{prp}
Let $\EE$ be a point-congruous $\Sigma$-Mal'tsev category. If $g.f$ and $g$ are $\Sigma$-special, so is $f:X\to Y$.
\end{prp}
\proof
The kernel congruence $R[f]$ is given by the following pullback in the category $Equ\EE$ of equivalence relations in $\EE$:
$$
\xymatrix@=15pt{
 {R[f]\;} \ar[rr]^{} \ar@{>->}[d]_{j} && {\Delta_{Y}\;}  \ar@{>->}[d]^{s_0}\\
 R[g.f] \ar[rr]_{R(f)} && R[g] 
             }
$$
where $\Delta_{Y}$ is the discrete equivalence relation on $Y$. The equivalence relations $R[g]$ and $R[g.f]$ are $\Sigma$-relations. Since the pullbacks in $Equ\EE$ are levelwise, and the class $\Sigma$ is point-congruous, the relation $R[f]$ is a $\Sigma$-relation as well.
\endproof
\begin{thm}
Let $\EE$ be a point-congruous $\Sigma$-Mal'tsev category. The subcategory $\Sigma l_Y\EE$ of the slice category $\EE/Y$ is a Mal'tsev category.
\end{thm}
\proof
Consider any reflexive relation $R$ in $\Sigma l_Y\EE$:
$$\xymatrix@=10pt
{
  R \ar[dd]_-{g_R} \ar@<2ex>[rr]^{d_{0}} \ar@<-2ex>[rr]_{d_{1}}& & X  \ar[dd]^-{g_X}  \ar[ll]_{s_0} \\
  &&& \\
  Y  \ar@{=}[rr] & & Y   
}
$$
According to the previous proposition the map $d_0$ is necessarily $\Sigma$-special. So according to Proposition \ref{sigrel}, the relation $R$ is an equivalence relation.
\endproof
In particular, if we denote by $\Sigma\EE_{\sharp}=\Sigma l_1\EE$ the full subcategory of $\EE$ whose objects are the $\Sigma$-special objects, it is a Mal'tsev category, called the \emph{Mal'tsev core} of the point-congruous $\Sigma$-Mal'tsev category $\EE$; any of its morphisms is $\Sigma$-special.

\begin{exm}\label{route}
\textbf{1)} The Mal'tsev core of the $\Sigma'$-Mal'tsev category $Mon$ of monoids is the category $Gp$ of groups, see \cite{BMMS}.\\
\textbf{2)} The Mal'tsev core of the $\bar{\Sigma}'$-Mal'tsev category $SRg$ of semi-rings is the category $Rg$ of rings, see \cite{BMMS}.\\
\textbf{3)} The Mal'tsev core of the $\Sigma'$-Mal'tsev category $Qnd$ of quandles is the category $LQd$ of latin quandles, see \cite{B15}.\\
\textbf{4)} The Mal'tsev core of the $\Sigma_Y$-Mal'tsev category $Cat_Y$ is the full subcategory of $Cat_Y$ whose objects are the categories $\mathbb Y$ having $Y$ as set of objects and such that, for all $y\in Y$, $Hom(y,y)$ is a group which acts in a simply transitive way on any non-empty $Hom(y,y')$, see \cite{B14}. This core is much larger than the one expected from example 1, namely the class of groupoids having $Y$ as set of objects.\\
\end{exm}

\section{The regular context}

In this section we shall suppose that the category $\EE$ is regular. Then the category $Pt(\EE)$ is regular as well, and its regular epimorphisms are the levelwise ones. We shall show that there are three meaningful levels of regular context for the $\Sigma$-Mal'tsev categories.

\subsection{$\P_{\EE}$-cartesian regular epimorphims}

First, we can extend Lemma \ref{transal} to any pullback of split epimorphims along a regular epimorphism:

\begin{lma}\label{transalreg}
Let $\EE$ be a regular $\Sigma$-Mal'tsev category. Consider any pullback of split epimorphism along a regular epimorphism $y$:
$$\xymatrix@=12pt
{
  X' \ar[dd]_{f'}  \ar@{->>}[rrrr]_{x} & &&& X  \ar[dd]_{f}  \\
  &&& \\
  Y' \ar@{->>}[rrrr]_{y}  \ar@<-1ex>[uu]_{s'} && && Y   \ar@<-1ex>[uu]_{s} 
}
$$
When $(f',s')$ is in $\Sigma$ the upward square is a pushout. Accordingly the base-change functor $y_{\Sigma}^*:\Sigma_Y(\EE)\to \Sigma_{Y'}(\EE)$ is fully faithful.
\end{lma}
\proof
The proof is exactly the same as the one of Lemma \ref{transal} since the regular epimorphism $x$ is again the quotient of its kernel relation.
\endproof

\begin{cor}
Let $\EE$ be a regular $\Sigma$-Mal'tsev category. Consider the following horizontal $\P_{\EE}$-cartesian regular epimorphism with domain $(f',s')\in Pt_{Y'}(\EE)$ in $\Sigma$:
$$\xymatrix@=12pt
{
 &&  \bar X \ar[rrd]^{x''} \ar[dddrr]_{\bar f}\\
  X' \ar[dd]_{f'}  \ar@{->>}[rrrr]_{x} \ar[rru]^{x'}& &&& X  \ar[dd]_<<<<{f}  \\
  &&& \\
  Y' \ar@{->>}[rrrr]_{y}  \ar@<-1ex>[uu]_{s'} && && Y   \ar@<-1ex>[uu]_>>>>{s} \ar@<-1ex>[uuull]_{}
}
$$
Suppose it has a decomposition through a map $x''$ in the fibre $Pt_Y\EE$. Then there is a splitting $m$ of $x''$ in this fibre such that $x'=m.x$. When, moreover, $x'$ is a regular epimorphism, then $x''$ is an isomorphism.
\end{cor}
\proof
The previous lemma says that the upward rectangle is a pushout; whence a unique factorization $m:X\into \bar X$ in $Pt_Y\EE$ such that $x'=m.x$. and $m.s=\bar s$. Then $x''.m.x=x''.x'=x$, and $x''.m=1_X$. When $x'$ is a regular epimorphism, so is the monomorphism $m$ which, accordingly, is an isomorphism.
\endproof

\begin{cor}
Let $\EE$ be a regular $\Sigma$-Mal'tsev category. Consider the following diagram where the whole rectangle is a pullback and its domain $(f',s')$ is in $\Sigma$:
$$\xymatrix@=4pt{
   {X'\;}   \ar@<-1ex>[ddd]_{f'}   \ar@{->>}[rrrrr]^{x} &&&&&{\bar X\;} \ar@<-1ex>[ddd]_{\bar f}\ar[rrrrr]^{\bar x} &&&&&  {X\;}  \ar@<-1ex>[ddd]_{f}\\
   &&&&&&&&\\
   &&&&&&&&\\
    Y' \ar@{->>}[rrrrr]_{y}  \ar[uuu]_{s'}  &&&&& \bar Y \ar[uuu]_{\bar s} \ar[rrrrr]_{\bar y}  &&&&& Y  \ar[uuu]_{s}
       }
       $$
       When $x$ is a regular epimorphisms, both squares are pullbacks.
\end{cor}
\proof
Set $(\check f,\check s)=\bar y^*(f,s)$, consider the following diagram, where $\check x$ is the canonical factorization induced by the whole rectangle and produces a pullback since the whole rectangle is itself a pullback:
$$\xymatrix@=12pt
{
 &&  \bar X \ar[rrd]^{x''} \ar[dddrr]_{\bar f}\\
  X' \ar[dd]_{f'}  \ar@{->>}[rrrr]_{\check x} \ar@{->>}[rru]^{x}& &&& \check X  \ar[dd]_<<<<{\check f}  \\
  &&& \\
  Y' \ar@{->>}[rrrr]_{y}  \ar@<-1ex>[uu]_{s'} && && \bar Y   \ar@<-1ex>[uu]_>>>>{\check s} \ar@<-1ex>[uuull]_{}
}
$$
Let $x''$ the factorization induced by the right hand side square. We have $x''.x=\check x$. According to the previous lemma, the map $x''$ is an isomorphism.
\endproof

\subsection{Regular pushouts}

We shall denote by $Reg^2(\EE)$ the category whose objects are regular epimorphisms and whose maps are the commutative squares between regular epimorphisms. In any regular category, it is straightforward that a pullback of a regular epimorphim along a regular epimorphism is a pushout.

\begin{dfn}
Let $\EE$ be a regular category. A commutative square of regular epimorphism as on the right hand side:
$$\xymatrix@=4pt{
{R[x]\;}   \ar@{->>}[ddd]_{R(f)}  \ar@<1ex>@{.>}[rrrrr]^{d_0^x} \ar@<-1ex>@{.>}[rrrrr]_{d_1^x} &&&&&{X\;}  \ar@{.>}[lllll] \ar@{->>}[ddd]_{f}  \ar@{->>}[rrrrr]^{x} &&&&&  {X'\;}  \ar@{->>}[ddd]^{f'}  \\
&&&&&&&&\\
&&&&&&&&\\
{R[y]\;}     \ar@<1ex>@{.>}[rrrrr]^{d_0^y} \ar@<-1ex>@{.>}[rrrrr]_{d_1^y} &&&&& Y \ar@{->>}[rrrrr]_{y} \ar@{.>}[lllll]  &&&&& Y' 
    }
    $$
is said to be a regular pushout whenever the induced factorization $X\to Y\times_{Y'} X'$ through the domain of the pullback of $f'$ along $y$ is a regular epimorphism.
\end{dfn}
\noindent Such a square is certainly a pushout and the factorization $R(f)$ is certainly a regular epimorphism.

\begin{lma}
Let $\EE$ be a regular category. The regular pushouts are stable under pullback along any map in $Reg^2(\EE)$.
\end{lma}
\proof
Consider a pair of commutative squares:
$$\xymatrix@=4pt{
 {X\;}   \ar@{->>}[ddd]_{f}   \ar@{->>}[rrrrr]^{x} &&&&&{X'\;} \ar@{->>}[ddd]^{f'} &&&& T \ar[ddd]_{g}  \ar@{->>}[rrrrr]^{t} &&&&&  {T'\;}  \ar[ddd]^{g'}\\
   &&&&&&&&\\
   &&&&&&&&\\
    Y \ar@{->>}[rrrrr]_{y}   &&&&& Y'  &&&& Y \ar@{->>}[rrrrr]_{y}   &&&&& Y'  
       }
       $$
    where the left hand side one (1) is a regular pushout and the right hand side one (2) is in $Reg^2(\EE)$. Set $\bar f=g^*(f):\bar X \onto T$ and $\bar f'=g'^*(f'):\bar X' \onto T'$. Consider now the following diagram:
    $$\xymatrix@=8pt{
    {\bar X\;} \ar[rrd]^{\bar{\phi}}  \ar@{->>}[ddd]_{\bar f}  \ar[rrrrr]^{\bar g} &&&&&  {X\;} \ar@{->>}[rrd]^{\phi}  \ar@{->>}[ddd]_>>>>>>{f}\\
    && T\times_{T'}\bar X' \ar@{->>}[ddll]^{t^*(\bar f')}  \ar[rrrrr]^<<<<<<<<<{}  &&&&& Y\times_{Y'}X'  \ar@{->>}[ddll]^{y^*(f')}  && \\
         &&& &&&&&&\\ 
    T \ar@{->>}[rrrrr]_{g}   &&&&& Y 
        }
        $$
      where the vertical square is a pullback by definition of $\bar f$ and the lower quadrangle as well by definition of $\bar f'$. So the upper quadrangle is a pullback. The comparison $\phi$ is a regular epimorphism since it comes from a regular pushout. Accordingly $\bar{\phi}$ is a regular epimorphism as well which means precisely that the pullback of (1) along (2) is a regular pushout.
\endproof

\begin{prp}\label{regpush}
Let $\EE$ be a regular $\Sigma$-Mal'tsev category. Consider a regular epimorphism in $Pt(\EE)$: 
$$\xymatrix@=4pt{
{X\;}   \ar@<-1ex>[ddd]_{f}  \ar@{->>}[rrrrr]^{x} &&&&&  {X'\;}  \ar@<-1ex>[ddd]_{f'}  \\
&&&&\\
&&&&\\
Y \ar@{->>}[rrrrr]_{y}  \ar[uuu]_{s}  &&&&& Y' \ar[uuu]_{s'} 
    }
    $$
Then, if its domain $(f,s)$ is in $\Sigma$, it is a regular pushout.
\end{prp}
\proof
We have to check that the factorization $\phi: X\to Y\times_{Y'}X'$ is a regular epimorphism. For that, complete the diagram with the kernel relations:
$$\xymatrix@=4pt{
{R[x]\;}   \ar@<-1ex>[ddd]_{R(f)}  \ar@<1ex>[rrrrr]^{d_0^x} \ar@<-1ex>[rrrrr]_{d_1^x} &&&&&{X\;}  \ar[lllll] \ar@<-1ex>[ddd]_{f}  \ar@{->>}[rrrrr]^{x} &&&&&  {X'\;}  \ar@<-1ex>[ddd]_{f'}  \\
&&&&&&&&\\
&&&&&&&&\\
{R[y]\;}   \ar[uuu]_{R(s)}  \ar@<1ex>[rrrrr]^{d_0^y} \ar@<-1ex>[rrrrr]_{d_1^y} &&&&& Y \ar@{->>}[rrrrr]_{y} \ar[lllll] \ar[uuu]_{s}  &&&&& Y' \ar[uuu]_{s'} 
    }
    $$
    Since $(f,s)$ is in $\Sigma$ the factorization $\bar{\phi}:R[x]\to R[y]\times_0X$, where $R[y]\times_0X$ is the domain of the pullback of $(f,s)$ along $d_0^y$, is a regular epimorphism. Since the following square commutes:
    $$\xymatrix@=15pt{
      R[x] \ar@{->>}[rr]^{d_1^x} \ar@{->>}[d]_{\bar{\phi}} && \ar[d]^{\phi}  X \\
      R[y]\times_0X \ar@{->>}[rr]_{d_1^y\times_0 x} && Y\times_{Y'} X
      }
      $$
    and since the left hand side and the lower maps are regular epimorphisms, certainly $\phi$ is a regular epimorphism as well.
    \endproof
    \noindent From that we can partially recover an important regular factorization result of the global Mal'tsev context, true for any split epimorphism $(f,s)$:
    \begin{cor}
   Let $\EE$ be a regular $\Sigma$-Mal'tsev category. Consider any pair of commutative squares with horizontal regular epimorphisms:
   $$\xymatrix@=4pt{
   {X\;}   \ar@<-1ex>[ddd]_{f}   \ar@{->>}[rrrrr]_{x} &&&&&{X'\;} \ar@<-1ex>[ddd]_{f'} &&&& T \ar[ddd]_{g}  \ar@{->>}[rrrrr]^{t} &&&&&  {T'\;}  \ar[ddd]^{g'}\\
   &&&&&&&&\\
   &&&&&&&&\\
    Y \ar@{->>}[rrrrr]_{y}  \ar[uuu]_{s}  &&&&& Y' \ar[uuu]_{s'} &&&& Y \ar@{->>}[rrrrr]_{y}   &&&&& Y'  
       }
       $$
      If $(f,s)$ is is $\Sigma$, then the factorization $x\times_yt:X\times_YT\to X'\times_{Y'}T'$ is a regular epimorphism.
    \end{cor}
    \proof
    Since $(f,s)$ is in $\Sigma$, the left hand side square above is a regular pushout. According to the previous lemma, its pullback along the right hand side square which is in $Reg^2(\EE)$ is a regular pushout:
         $$\xymatrix@=8pt{
               X\times_YT \ar[rr]^{x\times_yt} \ar@<-1ex>[dd]_{\bar{f}} &&  X'\times_{Y'}T'\ar@<-1ex>[dd]_{\bar{f}'} \\
               &&&\\
                T  \ar@{->>}[rr]_{t} \ar[uu]_{\bar s} &&  {T'\;}  \ar[uu]_{\bar s'}
               }
         $$
        and cerlainly its upper horizontal arrow is a regular epimorphism.
   \endproof
The direct image along a regular epimorphism of a reflexive (resp. symmetric) relation is reflexive (resp. symmetric). In general the direct image along a regular epimorphism of a transitive relation is no longer transitive, but this is the case in the Mal'tsev context. 
\begin{thm}
Let $\EE$ be a regular $\Sigma$-Mal'tsev category. The direct image along a regular epimorphim of a $\Sigma$-relation is transitive. A fortiori, the direct image along a regular epimorphim of a $\Sigma$-equivalence relation is an equivalence relation.
\end{thm}
\proof
Let $f:X\onto Y$ be a regular epimorphism. Consider the following double relation: 
$$ \xymatrix@=25pt{
      R[f] \square S \ar@<-1,ex>[d]_{p_0^f}\ar@<+1,ex>[d]^{p_1^f} \ar@<-1,ex>[r]_{p_0^S}\ar@<+1,ex>[r]^{p_1^S}
     & S \ar@<-1,ex>[d]_{d_0^S}\ar@<+1,ex>[d]^{d_1^S} \ar[l] \ar@{.>>}[r]^{\tilde f} & f(S) \ar@<-1,ex>[d]_{\delta_0}\ar@<+1,ex>[d]^{\delta_1}\\
       R[f] \ar@<-1,ex>[r]_{d_0^f} \ar@<+1,ex>[r]^{d_1^f} \ar[u]_{} & X \ar[u]_{} \ar[l] \ar@{->>}[r]_{f} & Y \ar[u]_{}
                     }
  $$
The upper equivalence relation is nothing but the kernel relation $R[f\times f.(d_0^S,d_1^S)]$, so it is an effective equivalence relation which has a quotient which is nothing but the direct image $f(S)$. When $S$ is a $\Sigma$-relation on $X$, on the one hand it is transitive and so is the left hand side vertical reflexive relation on $R[f]$, and on the other hand the right hand side square indexed by $0$ is a regular pushout by Proposition \ref{regpush}. Then complete the diagram by the pullback of the vertical $0$-indexed maps along the vertical $1$-indexed ones:
$$ \xymatrix@=25pt{
      (R[f] \square S)_2 \ar@<-2,ex>[d]_{p_0^f} \ar@{.>}[d]\ar@<+2,ex>[d]^{p_2^f} \ar@<-1,ex>[r]_{p_{02}}\ar@<+1,ex>[r]^{p_{12}}
           & S_2 \ar@<-2,ex>[d]_{d_0^S}\ar@{.>}[d]\ar@<+2,ex>[d]^{d_1^S} \ar[l] \ar@{->>}[r]^{\tilde f_2} & f(S)_2 \ar@<-2,ex>[d]_{\delta_0}\ar@{.>}[d]\ar@<+2,ex>[d]^{\delta_1}\\
      R[f] \square S \ar@<-1,ex>[d]_{p_0^f}\ar@<+1,ex>[d]^{p_1^f} \ar@<-1,ex>[r]_{p_0^S}\ar@<+1,ex>[r]^{p_1^S}
     & S \ar@<-1,ex>[d]_{d_0^S}\ar@<+1,ex>[d]^{d_1^S} \ar[l] \ar@{->>}[r]^{\tilde f} & f(S) \ar@<-1,ex>[d]_{\delta_0}\ar@<+1,ex>[d]^{\delta_1}\\
       R[f] \ar@<-1,ex>[r]_{d_0^f} \ar@<+1,ex>[r]^{d_1^f} \ar[u]_{} & X \ar[u]_{} \ar[l] \ar@{->>}[r]_{f} & Y \ar[u]_{}
                     }
  $$
 By commutations of limits, the upper row produces the kernel relation of $\tilde f_2$. The main fact is that, according to the previous corollary, the factorization $\tilde f_2$ is a regular epimorphism which is consequently the quotient of this kernel relation. This allows to transfer the vertical dotted ``transitivity'' arrows of $S$ and of the vertical reflexive relation on $R[f]$ to the ``quotient'' relation $f(S)$.
\endproof

\subsection{Three levels of regular context}

The $\Sigma$-Mal'tsev categories have three possible degrees of intricateness with regularness which are exemplified by different contexts.

Consider any regular epimorphism in $Pt(\EE)$ with its domain $(f,s)$ in $\Sigma$:
$$\xymatrix@=4pt{
{X\;}   \ar@<-1ex>[ddd]_{f}  \ar@{->>}[rrrrr]^{x} &&&&&  {X'\;}  \ar@<-1ex>[ddd]_{f'}  \\
&&&&\\
&&&&\\
Y \ar@{->>}[rrrrr]_{y}  \ar[uuu]_{s}  &&&&& Y' \ar[uuu]_{s'} 
    }
    $$
\begin{dfn}
Let $\EE$ be a regular $\Sigma$-Mal'tsev category. We shall say it is:\\
1) 1-regular when $(f',s')$ is in $\Sigma$ for any $\P_{\EE}$-cartesian regular epimorphism as above\\
2) 2-regular when $(f',s')$ is in $\Sigma$ for any regular epimorphism as above when the split epimorphism $(R(f),R(s))$ induced on the kernel relations is in $\Sigma$ as well\\
3) 3-regular when $(f',s')$ is in $\Sigma$ for any regular epimorphism as above.
\end{dfn}

\noindent It is clear that these notions are of increasing strength.

\begin{exm}
1) Proposition 2.3.5 in \cite{BMMS} shows that the category $Mon$ of monoids is 2-regular with respect to the class $\Sigma'$\\
2) it is straighforward to check that the category $Mon$ of monoids is 3-regular with respect to the class $\Sigma$\\
3) suppose that $U:\CC\to \DD$ is a left exact conservative functor and $\Sigma$ is a fibrational class of split epimorphisms. Set $\bar{\Sigma}=U^{-1}(\Sigma)$. Suppose moreover that $U$ preserves and reflects regular epimorphisms. Then $\CC$ is i)-regular with respect to $\bar{\Sigma}$ as soon as so is $\DD$ with respect to $\Sigma$\\
4) it is the case for the forgetful functor $U:SRg\to Mon$; so that the category $SRg$ of semi-rings is 2-regular with respect to the class $\bar{\Sigma}'$ and 3-regular with respect to the class $\bar{\Sigma}$\\
5) the category $Qnd$ of quandles is 3-regular with respect to the class $\Sigma$ and 1-regular with respect to $\Sigma'$, see Proposition 2.6 in \cite{B15}\\
6) let $FQnd$ denotes the category of finite quandles; the classes $\Sigma$ and $\Sigma'$ coincide in $FQnd$, so that this category is, at the same time, a 3-regular and a point-congruous $\Sigma$-Mal'tsev category.	
\end{exm}

\begin{prp}
Let $\EE$ be a regular $\Sigma$-Mal'tsev category which is 1-regular. Then pulling back along regular epimorphims reflects the split epimorphisms in $\Sigma$ and reflects the $\Sigma$-special morphisms.
\end{prp}
\proof
The level 1 of regularness is exactly the first part of the assertion. To get the second one, extend the pullback in question by its kernel relations.
\endproof

\begin{prp}
Let $\EE$ be a regular point-congruous $\Sigma$-Mal'tsev category which is 2-regular. Then the category $\Sigma(\EE)$ is a regular category whose regular epimorphisms are the levelwise ones. The same result holds for any Mal'tsev fibre $\Sigma l_Y(\EE)$.
\end{prp}
\proof
Consider any morphism in $\Sigma(\EE)$ as on the left hand side:
$$
\xymatrix@=20pt{
R[x] \ar@<-1,ex>@{.>}[rr]_{d_0^x} \ar@<+1,ex>@{.>}[rr]^{d_1^x} \ar@<-4pt>[d]_{R(\bar f)} && \bar X  \ar[rr]^{x}\ar@{.>}[ll] \ar@<-4pt>[d]_{\bar f} && X \ar@<-4pt>[d]_{f}  && \bar X  \ar@{->>}[rr]^{\epsilon_x} \ar@<-4pt>[d]_{\bar f} && {X'\;} \ar@{>->}[rr]^{m_x} \ar@<-4pt>[d]_{f'} && X \ar@<-4pt>[d]_{f} \\
R[y] \ar@<-1,ex>@{.>}[rr]_{d_0^y} \ar@<+1,ex>@{.>}[rr]^{d_1^y} \ar[u]_{R(\bar s)} && \bar Y \ar[rr]_{y} \ar@{.>}[ll] \ar[u]_{\bar s} && Y \ar[u]_{s} &&   \bar Y \ar@{->>}[rr]_{\epsilon_y}  \ar[u]_{\bar s} && {Y'\;} \ar@{>->}[rr]_{m_y} \ar[u]_{s'} && Y \ar[u]_{s}
}
$$
Complete it by its kernel relation. Since $\Sigma$ is point-congruous, the split epimorphism $(R(\bar f),R(\bar s))$ is in $\Sigma$. So the regular decomposition in $\EE$ on the right hand side is such that $(f',s')$ is in $\Sigma$ when $\EE$ is 2-regular. Accordingly the morphism in $\Sigma(\EE)$ we started with is regular if and only if both $y$ and $x$ are regular epimorphisms. Such morphisms are stable under pullbacks in $\Sigma(\EE)$ since $\EE$ is regular and $\Sigma$ is point-congruous. Now, the first part of the proof shows that any effective equivalence relation in $\Sigma(\EE)$ has a quotient.

Let $h$ be any morphism in the fibre $\Sigma l_Y(\EE)$ and consider its regular decomposition in the category $\EE$:
$$
\xymatrix@=25pt{
{X\;} \ar@<2ex>[rr]^{h} \ar@{->>}[r]_{\bar h}\ar[d]_{f} & {U\;} \ar@{>->}[r]_{m} \ar[d]^{f'.m} & X'  \ar[d]^{f'}\\
Y \ar@{=}[r] & Y \ar@{=}[r] & Y 
              }
$$
Then according to the previous part of the proof, if $R[f]$ and $R[f']$ are $\Sigma$-relations, so is $R[f'.m]$ and $f'.m$ is $\Sigma$-special. Accordingly the regular epimorphisms in the fibre $\Sigma l_Y(\EE)$ are those morphisms which are given by a morphism $h$ which is a regular epimorphism in $\EE$. The end of the proof is then exactly as above.
\endproof

\subsection{Baer sums}

We shall show in this section that the level 1 of regularness allows to extend to the partial Mal'tsev context the construction of the Baer sum of extensions with abelian kernel relation valid in the global regular Mal'tsev context, provided that we restrict our attention to the $\Sigma$-special extensions. This will include the Baer sum of special Schreier extensions with abelian kernel in Monoids and the Baer sum of special Schreier extensions with trivial kernel in Semi-rings as described in \cite{BMMS}.

More precisely we shall suppose in this section that \emph{the category $\EE$ is an efficiently regular $\Sigma$-Mal'tsev category which is regular of level 1}. We shall call \emph{abelian $\Sigma$-special extension} any $\Sigma$-special regular epimorphism $f$ with abelian kernel relation, namely such that $[R[f],R[f]]=0$; in this case consider the following diagram given by the double centralizing relation:
$$ \xymatrix@=30pt{
      R[f]\rtimes_XR[f] \ar@<-1,ex>[d]_{\pi_0^f}\ar@<+1,ex>[d]^{\pi_1^f} \ar@<-1,ex>[r]_{p_0^f}\ar@<+1,ex>[r]^{p_1^f}
     & R[f] \ar@<-1,ex>[d]_{d_0^f}\ar@<+1,ex>[d]^{d_1^f} \ar[l] \ar@{.>>}[r]^{q_f} & \bar A \ar@{.>>}@<+1,ex>[d]^{\bar f}\\
       R[f] \ar@<-1,ex>[r]_{d_0^f} \ar@<+1,ex>[r]^{d_1^f} \ar[u]_{} & X\ar@{->>}[r]_{f}\ar[l]\ar[u]_{}  & Y \ar@{.>}[u]^{\bar s} 
                     }
  $$
Since $\EE$ is efficiently regular and since the left hand side square indexed by $0$ is a discrete fibration between equivalence relation, the upper equivalence relation is effective, and so has a quotient $q_f$ which produces the split epimorphism $(\bar f,\bar s)$ and makes the right hand side square a pullback.
\begin{prp}
The split epimorphism $(\bar f,\bar s)$ is $\Sigma$-special and has an abelian kernel relation. It is called the \emph{direction} of the $\Sigma$-special extension $f$.
\end{prp}
\proof
The split epimorphism $(\bar f,\bar s)$ is $\Sigma$-special since the right hand side square is a pullback and pulling back along a regular epimorphism reflects the $\Sigma$-special morphisms. Moreover the vertical right hand side part is necessarily a group in the fibre $\Sigma l_Y(\EE)$ as a quotient of the vertical groupoid (actually an equivalence relation) $R[f]$. This group is necessarily abelian by Corollary \ref{abgp}. Accordingly $\bar f$ has an abelian kernel relation.
\endproof
Suppose now given any abelian group $(\bar f,\bar s):\bar A\rightleftarrows Y$ in $\Sigma _Y(\EE)$. For sake of simplicity we shall denote this whole abelian structure by the only symbol $\bar A$. We shall call $\bar A$-torsor a triple $(f,R[f],q_f)$ of a regular epimorphism $f$ together with a (regular) discrete fibration between the following vertical internal groupoids:
$$ \xymatrix@=30pt{
     R[f] \ar@<-1,ex>[d]_{d_0^f}\ar@<+1,ex>[d]^{d_1^f} \ar@{->>}[r]^{q_f} & \bar A \ar@{->>}@<+1,ex>[d]^{\bar f}\\
      X\ar@{->>}[r]_{f}\ar[u]_{}  & Y \ar[u]^{\bar s} 
                     }
  $$
The split epimorphism $(\bar f,\bar s)$ being in $\Sigma$ and the previous square being a pullback, the regular epimorphism $f$ is $\Sigma$-special and we do get $[R[f],R[f]]=0$. Conversely, according to the previous construction, any abelian $\Sigma$-special extension gives rise to a $\bar A$-torsor. A morphism of $\bar A$-torsors $(f,R[f],q_f)\to (f',R[f'],q_{f'})$ is given by a map $h:X\to X'$ such that $f'.h=f$ and $q_{f'}.R(h)=q_f$.
Any $\bar A$-torsor $(f,R[f],q_f)$ produces another one, namely $(f,R^{op}[f],q_f)$, where the projection are twisted:
$$ \xymatrix@=30pt{
     R[f] \ar@<-1,ex>[d]_{d_1^f}\ar@<+1,ex>[d]^{d_0^f} \ar@{->>}[r]^{q_f} & \bar A \ar@{->>}@<+1,ex>[d]^{\bar f}\\
      X\ar@{->>}[r]_{f}\ar[u]_{}  & Y \ar[u]^{\bar s} 
                     }
  $$
In set theoretical terms, the groupoid structure given on $R[f]$ implies that $q_f(a,b)+q_f(b,a)=q_f(a,a)=0$, in other words, that $q_f(b,a)=-q_f(a,b)$, or equivalently that $q_f.tw=-1_{\bar A}.q_f$, where $tw:R[f]\to R^{op}[f]$ is the twisting isomorphism. Accordingly, the homomorphism induced on $\bar A$ by the twisting isomorphism $tw$ is nothing but $-1_{\bar A}$.

Let us denote by $EXT(\bar A,Y)$ the category whose objects are the $\bar A$-torsors and whose morphisms $h$ are the morphism of $\bar A$-torsors which induce $1_{\bar A}$ on $\bar A$. As usual we get:
\begin{cor}
Any morphism of $\bar A$-torsors is an isomorphism.
\end{cor}
\proof
Consider the following diagram where $h$ is a morphism of $\bar A$-torsors: 
  $$\xymatrix@=8pt{
{R[f]\;} \ar[rrd]^{R(h)}  \ar@<-1ex>[ddd]_{d_0^f} \ar@<1ex>[ddd]^{d_1^f}  \ar[rrrrr]^{q_f} &&&&&  {\bar A\;} \ar@{=}[rrd]   \ar@<1ex>[ddd]^{\bar f}\\
&& R[f'] \ar@<-1ex>[ddd]_<<<<{d_0^{f'}} \ar@<1ex>[ddd]^<<<<{d_1^{f'}} \ar[rrrrr]^<<<<<<<<<{q_{f'}}  &&&&& \bar A' \ar@<1ex>[ddd]^{\bar f'}  && \\
     &&& &&&&&&\\ 
X \ar[rrd]_{h}  \ar@{->>}[rrrrr]^>>>>>{f} \ar@{->>}[dd]_f \ar[uuu]_{}  &&&&& Y \ar[uuu]_{} \ar@{=}[rrd]  \\
&&  X'  \ar@{->>}[rrrrr]_{f'} \ar@{->>}[dd]_{f'}  \ar[uuu]_{}  &&&&& Y \ar[uuu]_{} \\
Y \ar@{=}[rrd] && \\
&& Y'
    }
    $$
    Then the left hand side quadrangle is a pullback. By the Barr-Kock theorem, the lower square is a pulback as well, and $h$ is an isomorphism.
\endproof
Accordingly the category $EXT(\bar A,Y)$ is a groupoid. Let us denote by $Ext(\bar A,Y)$ of set of connected components of this groupoid.
Classically there is a symmetric tensor product on the $\bar A$-torsor, see for instance \cite{Ba}.  This tensor product allows us to define the \emph{Baer sum} on the set $Ext(\bar A,Y)$ which gives it an abelian group structure. Starting with two $\bar A$-torsors $f$ and $f'$ in $Ext(\bar A,Y)$, the Baer sum is given by the following construction; take the following left hand side quadrangled pullbacks, paying attention to the fact that the upper relation is $R[f']$ while the lower one is $R^{op}[f]$:
  $$\xymatrix@=8pt{
{R\;}\ar[ddd]^{} \ar[rrd]^{}  \ar@<-1ex>[rrrrr]_{d_0} \ar@<1ex>[rrrrr]^{d_1}   &&&&&  {X\times_YX'\;} \ar[rrd]^{p_{X'}}  \ar[lllll] \ar[ddd]^{} \ar@{.>>}[rrrrr]^{q} &&&&& X\otimes X' \ar@{->>}[drr]^{f\otimes f'} \ar@{->>}[ddd]^>>>>>{f\otimes f'}\\
&& R[f'] \ar@<-1ex>[rrrrr]_<<<<<<<<{d_0^{f'}} \ar@<1ex>[rrrrr]^<<<<<<<<{d_1^{f'}} \ar[ddd]_<<<<<{q_{f'}} &&&&& X' \ar@{->>}[rrrrr]^<<<<<<<<<{f'}  \ar@{->>}[ddd]_<<<<<{f'} \ar[lllll] &&&&& Y \\
     &&& &&&&&& &&&\\ 
{R[f]\;} \ar[rrd]_{q_f}  \ar@<-1ex>[rrrrr]_>>>>>{d_1^f} \ar@<1ex>[rrrrr]^>>>>>{d_0^f}   &&&&& X\ar[lllll]  \ar@{->>}[rrrrr]_>>>>>{f} \ar@{->>}[rrd]^f  &&&&& Y\\
&&  \bar A  \ar@<-1ex>[rrrrr]_{\bar f}   &&&&& Y \ar[lllll]_{\bar s} &&&&\\
}
    $$
These pullbacks define an equivalence relation $R$ on $X\times_YX'$. Since the back right hand side part of the diagram is underlying a discrete fibrations between equivalence relations and the category $\EE$ is efficiently regular, the equivalence relation $R$ has a quotient $q$ and produces a unique factorization $f\otimes f'$ such that $f\otimes f'.q=p_X.f=p_{X'}.f'$. It is a $\Sigma$-special morphism since, in a regular $\Sigma$-Mal'tsev category of level 1, pulling back reflects this kind of morphisms. This is the Baer sum of $f$ and $f'$. As usual, the inverse of $(f,R[f],q_f)$ is precisely $(f,R^{op}[f],q_f)$.
 
\section{Partial protomodularity}

\subsection{Strongly split epimorphism}

\noindent Let $\EE$ be a finitely complete category.

\begin{dfn}
A split epimorphism $(f,s)$ in $\EE$ is called strongly split, when, given any pullback:
$$
\xymatrix@=20pt{
\bar X  \ar[rr]^{x} \ar@<-4pt>[d]_{\bar f} && X \ar@<-4pt>[d]_{f}  && \\
\bar Y \ar[rr]_{y}  \ar[u]_{\bar s}& & Y \ar[u]_{s}}
$$
the pair $(x,s)$ is jointly extremally epic.
\end{dfn}

When the category $\EE$ is pointed a split epimorphism is strongly split if and only if in the following pullback:
$$
\xymatrix@=20pt{
{K[f]\;}  \ar@{>->}[rr]^{k_f} \ar@<-4pt>[d] && X \ar@<-4pt>[d]_{f}  && \\
{1\;} \ar@{>->}[rr]_{\alpha_Y}  \ar[u] & & Y \ar[u]_{s}}
$$
the pair $(k_f,s)$ is jointly extremally epic (or, equivalently the subobject $1_X$ is the supremum of the subobjects $s$ and $k_f$). Accordingly, in the pointed context, the notion of strongly split epimorphism introduced here is equivalent to the one defined in \cite{BMMS}.

When $\EE$ is a regular category (which includes all the varietal examples), we have again the following description in terms of supremum:

\begin{prp}
Let $\EE$ be a regular category. The split epimorphism $(f,s)$ is strongly split if and only if, for any pullback above a monomorphism $m$:
$$
\xymatrix@=20pt{
{X'\;} \ar@{>->}[rr]^{n} \ar@<-4pt>[d]_{f'} && X \ar@<-4pt>[d]_{f}  \\
{Y'\;} \ar@{>->}[rr]_{m}  \ar[u]_{s'}& & Y \ar[u]_{s}}
$$
the pair $(s,n)$ of monomorphisms is jointly strongly epic, or in other words, if and only if the subobject $1_X$ is the supremum of the subobjects $s$ and $n$.
\end{prp}
\proof
Consider any pullback on the left hand side and its canonical regular decomposition on the right hand side:
$$
\xymatrix@=20pt{
\bar X  \ar[rr]^{x} \ar@<-4pt>[d]_{\bar f} && X \ar@<-4pt>[d]_{f}  && \bar X  \ar@{->>}[rr]^{\epsilon_x} \ar@<-4pt>[d]_{\bar f} && {X'\;} \ar@{>->}[rr]^{m_x} \ar@<-4pt>[d]_{f'} && X \ar@<-4pt>[d]_{f} \\
\bar Y \ar[rr]_{y}  \ar[u]_{\bar s} && Y \ar[u]_{s} &&   \bar Y \ar@{->>}[rr]_{\epsilon_y}  \ar[u]_{\bar s} && {Y'\;} \ar@{>->}[rr]_{m_y} \ar[u]_{s'} && Y \ar[u]_{s}
}
$$
Since $\EE$ is regular, the two right hand squares are pullback as well. Now, the pair $(x,s)=(m_x.\epsilon_x,s)$ is jointly strongly epic if and only if so is the pair $(m_x,s)$ since $\epsilon_x$ is a regular epimorphism. It is the case by assumption since the square above $m_y$ is a pullback.
\endproof

\subsection{$\Sigma$-protomodularity}

Recall that a category $\EE$ is said to be \emph{protomodular} when any base-change functor with respect to the fibration of points is conservative \cite{B}. It is equivalent to say that any split epimorphism is strongly split. Recall also that any protomodular category is a Mal'tsev one \cite{B0}.
\begin{dfn}\label{sproto}
Let $\EE$ be a category endowed with a fibrational class $\Sigma$. It is be said to be $\Sigma$-protomodular when any split epimorphism in $\Sigma$ is strongly split.
\end{dfn}
\begin{thm}\label{main}
   Let $\EE$ be a category endowed with a fibrational class $\Sigma$ of split epimorphisms. Then:\\
   1) when $\EE$ is $\Sigma$-protomodular, it is a $\Sigma$-Mal'tsev category,\\
   2) when, in addition, $\Sigma$ is point-congruous, any base-change functor with respect to the fibration $\P_{\EE}^{\Sigma}$ is conservative.
  \end{thm}
  \proof
  1) Consider the following (rightward) pullback of split epimorphisms in $\EE$ where $(f,s)$ in $\Sigma$:
  $$
    \xymatrix{
     X'  \ar@<-4pt>[r]_(.6){g'} \ar@<-4pt>[d]_{f'} & X \ar@<-4pt>[l]_(.4){t'} \ar@<-4pt>[d]_f\\
     Y'   \ar@<-4pt>[r]_g \ar@<-4pt>[u]_{s'} & Y \ar@<-4pt>[l]_t \ar@<-4pt>[u]_s }
  $$
  The split epimorphism $(f',s')$ is then in $\Sigma$ and so a strongly split epimorphism. Now, because of the splitting $t$, the leftward square is still a pullback of split epimorphisms, and the pair $(t',s)$ is certainly joinly strongly epic.
  
   2) When $\Sigma$ is point-congruous, any fibre $\Sigma_Y(\EE)$ is stable under finite limit in $Pt_Y\EE$. Then any base-change with respect to the fibration $\P_{\EE}^{\Sigma}$ is left exact, and it is enough to prove that it is conservative on monomorphisms. So let us consider the following diagram where all the quadrangles are pullbacks and all the split epimorphisms are in $\Sigma(\EE)$:
    $$\xymatrix@=12pt{
      {X'\;} \ar[rd]^{m'}_{\simeq}  \ar@<-1ex>[ddd]_{f'}  \ar[rrrr]^{x}
  &&&&  {X\;} \ar@{>->}[rd]^{m}  \ar@<-1ex>[ddd]_{f}  \\
  & \bar X'
  \ar[ddl]_{\bar f'} \ar[rrrr]^<<<<{\bar x}  & & &  & \bar X
  \ar[ddl]_{\bar f}  && \\
       &&& &&&&&&\\
   Y' \ar[rrrr]_{y}  \ar[uuu]_{s'}
  \ar@<-1ex>[uur]_{\bar s'} && && Y \ar[uuu]_{s}
  \ar@<-1ex>[uur]_{\bar s}
      }
      $$
      Suppose moreover that the factorization $m'$ is an isomorphism. Since the split epimorphim $(\bar f,\bar s)$ is a strongly split epimorphism the pair $(\bar x,\bar s)$ is jointly extremally epic; accordingly, since $m'$ is an isomorphism, so is the pair $(\bar x.m',\bar s)$, and $m$ is necessarily an isomorphism.
    \endproof
    
    \subsection{Examples}
 According to our previous remark on strongly split epimorphisms in a pointed context, the notion of $\Sigma$-protomodular category introduced for the pointed and point-congruous context of the category $Mon$ in \cite{BMMS}  coincides with the present one provided that the class $\Sigma$ is point-congruous. Whence the following two first examples of partial protomodularity:\\
1) the category $Mon$ of monoids is $\Sigma'$ protomodular\\
2) the category $SRg$ of semi-rings is $\bar{\Sigma}'$-protomodular\\
3) suppose that $U:\CC\to \DD$ is a left exact conservative functor and $\Sigma$ is a fibrational class of split epimorphisms. Set $\bar{\Sigma}=U^{-1}(\Sigma)$. Then $\CC$ is $\bar{\Sigma}$-protomodular as soon as $\DD$ is  $\Sigma$-protomodular\\
4) it is straightforward to check that the category $Mon$ is actually $\Sigma$-protomodular and the category $SRg$ is $\bar{\Sigma}$-protomodular\\
5) the fibres $Cat_Y(\EE)$ are non-pointed $\Sigma_Y$-protomodular categories, see \cite{B14}\\
6) because of the presence of $\emptyset$ as an object in the category $Qnd$ of quandles, this category is a discriminating example of a $\Sigma$-Mal'tsev category which is not $\Sigma$-protomodular.

\subsection{Properties}

In a protomodular category, there is an induced notion of normal monomorphism; here we have similarly:
\begin{prp}\label{normal}
Let $\EE$ be a $\Sigma$-Mal'tsev category. When $m:U\rightarrowtail X$ is a monomorphism which is normal to a $\Sigma$-equivalence relation $S$ on $X$, the object $U$ is $\Sigma$-special. When $\EE$ is $\Sigma$-protomodular, a monomorphism $m$ is normal to at most one $\Sigma$-equivalence relation. 
\end{prp}
\proof
Saying that the monomorphism $m$ is normal to an equivalence relation $S$ is saying that $m^{-1}(S)=\nabla_U$ (the indiscrete relation on $U$) and that the induced following functor is a discrete fibration, i.e. that any of the following commutative squares is a pullback:
$$  \xymatrix@=4pt{
{U\times U\;}   \ar@<-1ex>[ddd]_{p_0} \ar@<1ex>[ddd]^{p_1}  \ar@{>->}[rrrrr]^{\tilde m} &&&&&  {S\;}  \ar@<-1ex>[ddd]_{d_0}  \ar@<1ex>[ddd]^{d_1}\\
&&&&\\
&&&&\\
 {U\;} \ar@{>->}[rrrrr]_{m}  \ar[uuu]_{}  &&&&& X \ar[uuu]_{} 
    }
$$
In set theoretical terms, it is saying that $U$, when it is non-empty, is an equivalence class of $R$.
Accordingly, when $S$ is a $\Sigma$-equivalence relation, so is $\nabla_U$, and $U$ is a $\Sigma$-special object.
 
When a monomorphism $m$ is normal to two equivalence relations $R$ and $R'$, it normal to $R\cap R'$. Now suppose that $\EE$ is $\Sigma$-protomodular and consider the following diagram:
$$\xymatrix@=12pt{
  {U\times U\;} \ar@{=}[rd]  \ar@<-2ex>[ddd]_{p_0}  \ar@{>->}[rrrr]_{}  &&&&  {R\cap S\;} \ar@{>->}[rd]^j  \ar@<-1ex>[ddd]_{d_0}  \\
   & {U\times U\;} \ar[ddl]_{p_0} \ar@{>->}[rrrr]^<<<<{\tilde m}  & & &  & S \ar[ddl]_{d_0}  && \\
    &&&&&&\\
  {U\;} \ar@{>->}[rrrr]_{m}  \ar@<1ex>[uuu]_{s_0}  \ar@<-1ex>[uur]_{s_0} && && X \ar[uuu]_{s_0} \ar@<-1ex>[uur]_{s_0}
  }
  $$
  When $S$ is a $\Sigma$-equivalence relation, the downward quadrangle being a pullback, the pair $(\tilde m, s_0)$ is extremally epic, so that the monomorphism $j$ is an isomorphism, and we get $S\subset R$. If $R$ is a $\Sigma$-equivalence relation as well, we get $R=S$.
 \endproof

\noindent When, in addition, $\Sigma$ is point-congruous, we get:

\begin{prp}
Let $\EE$ be a category endowed with a point-congruous class $\Sigma$ and suppose it is $\Sigma$-protomodular. Then any fibre $\Sigma l_Z(\EE)$ is protomodular. In particular, the full subcategory $\Sigma\EE_{\sharp}=\Sigma l_1\EE$ of $\EE$ is protomodular; we call it the protomodular core of $\EE$.
\end{prp}
\proof
Consider the following diagram in $\Sigma l_Z(\EE)$:
$$\xymatrix@=12pt{
      {X'\;} \ar[rd]^{m'}_{\simeq}  \ar@<-1ex>[ddd]_{f'}  \ar[rrrr]^{x}
  &&&&  {X\;} \ar@{>->}[rd]^{m}  \ar@<-1ex>[ddd]_{f}  \\
  & \bar X' \ar[ddl]_{\bar f'} \ar[rrrr]^<<<<{\bar x}  & & &  & \bar X \ar[ddl]_{\bar f}  && \\
       &&& &&&&&&\\
   Y' \ar[rrrr]_{y}  \ar[uuu]_{s'} \ar@<-1ex>[uur]_{\bar s'}\ar[d]_{g'} && && Y \ar[uuu]_{s} \ar@<-1ex>[uur]_{\bar s} \ar[d]_{g}\\
   Z \ar@{=}[rrrr] &&&& Z
      }
      $$
Since any of the split epimorphisms is in $\Sigma l_Z(\EE)$, they are $\Sigma$-special split epimorphisms; so they are in $\Sigma$ thanks to Lemma \ref{splitsp}.
On the other hand since the base-change $y^*$ is left exact, it is enough to prove that it is conservative on monomorphisms. Then we can apply the point 2) of Theorem \ref{main}.
      \endproof
So, all the examples of core given in \ref{route} are protomodular, except the third one.

\noindent keywords: Mal'tsev and protomodular categories, fibration of points, permutation and centralization of equivalence relations, Baer sum.

\smallskip

\noindent Mathematics Subject Classification: 18C99, 18A10, 08B05, 57M27.

\smallskip

\noindent{\small\sc  Universit\'e du Littoral, Lab. de Math\'ematiques Pures et Appliqu\'ees, Bat. H. Poincar\'e, 50 Rue F. Buisson, BP 699, 62228 Calais Cedex, France}\\ \emph{E-mail:} bourn@lmpa.univ-littoral.fr

\end{document}